\documentclass[12pt]{amsart}

\usepackage{amsbsy,amsthm,amsfonts,amsmath,amssymb,eucal,amscd,graphics,enumerate, stmaryrd,xspace,verbatim, epic, eepic,epsfig,verbatim}
\usepackage[all]{xy}

\newcommand{\iX}[1]{{\mathcal I}_{\mu_r}({\mathcal X})}

\def\clA{{\mathcal A}}  \def\clC{{\mathcal C}} \def\clD{{\mathcal D}} 
 \def\clG{{\mathcal G}}  \def\clI{{\mathcal I}} \def\clK{{\mathcal K}}
 \def\clL{{\mathcal L}}
\def\clM{{\mathcal M}} \def\clN{{\mathcal N}} \def\clO{{\mathcal O}} \def\clP{{\mathcal P}}

  \def\clU{{\mathcal U}} 
  
\def\clX{{\mathcal X}}
\def\clY{{\mathcal Y}}

  \def\bbC{{\mathbb C}}

\def\bbQ{{\mathbb Q}} \def\bbZ{{\mathbb Z}}   \def\bbN{{\mathbb N}}

\def\ttau{{\tilde{\tau}}}

\def\bix{{\overline{I}_{\mu}(\clX)}}
\def\biy{{\overline{I}_{\mu}(\clY)}}
\def\bix1{{\overline{I}_{\mu_{r_1}}(\clX)}}
\def\biy2{{\overline{I}_{\mu_{r_2}}(\clY)}}

   \def\frM{\mathfrak{M}}
 
 \def\frB{\mathfrak{B}}  \def\frD{\mathfrak{D}}

\newtheorem{definition}{Definition}[section]

\newtheorem{prop}[definition]{Proposition}

\newtheorem{lemm}[definition]{Lemma}

\newtheorem{theo}[definition]{Theorem}
\newtheorem{remark}[definition]{Remark}

\newtheorem{notation}[definition]{Notation}

\newtheorem{remarkdef}[definition]{Remark-Definition}

\newenvironment{rema}{\begin{remark} \rm}{\end{remark}}

\newtheorem{remarks}[definition]{Remarks}

\newtheorem{example}[definition]{Example}

\newtheorem{examples}[definition]{Examples}

\newtheorem{nothing}[definition]{$\!\!$}

\newenvironment{proo}{{\flushleft \bf Proof :}}{\hfill $\square$ \vspace{5mm}}

\newtheorem{definition*}{Definition}[section]
\newenvironment{defi*}{\begin{definition*} \rm}{\end{definition*}}
\newtheorem{prop*}[definition*]{Proposition}
\newtheorem{lemm*}[definition*]{Lemma}
\newtheorem{coro*}[definition*]{Corollary}
\newtheorem{theo*}[definition*]{Theorem}
\newtheorem{remark*}[definition*]{Remark}
\newenvironment{rema*}{\begin{remark*} \rm}{\end{remark*}}
\newtheorem{remarks*}[definition*]{Remarks}
\newenvironment{remas*}{\begin{remarks*} \rm}{\end{remarks*}}
\newtheorem{example*}[definition*]{Example}
\newenvironment{exam*}{\begin{example*} \rm}{\end{example*}}
\newtheorem{examples*}[definition*]{Examples}\begin{large}                                             \end{large}
\newenvironment{exams*}{\begin{examples*} \rm}{\end{examples*}}

\def\longto{\longrightarrow}
\def\isomto{\stackrel{\sim}{\longto}}

\renewcommand{\tilde}{\widetilde}

\newcommand{\Spec}{\operatorname{Spec}}

\providecommand{\abs}[1]{\lvert#1\rvert}

\title{Gromov-Witten theory of product stacks}
\date{\today}
\author{Elena Andreini}
\address{Max-Planck-Institut f\"ur Mathematik\\Vivatsgasse 7\\
53111 Bonn\\ Germany}
\email{andreini.elena@gmail.com}
\author{Yunfeng Jiang}
\address{Department of Mathematics\\ University of Utah\\ 155 S 1400 E JWB 233\\Salt Lake city\\ UT 84112\\ USA}
\email{jiangyf@math.utah.edu}
\author{Hsian-Hua Tseng}
\address{Department of Mathematics\\ University of Wisconsin-Madison\\ Van Vleck Hall, 480 Lincoln Drive \\Madison\\ WI 53706-1388 \\ USA}
\email{tseng@math.wisc.edu}

\begin{document}
\sloppy \maketitle


\begin{abstract}
Let $\mathcal{X}_1$ and $\mathcal{X}_2$ be smooth proper Deligne-Mumford stacks with projective coarse moduli spaces. We prove a formula for orbifold Gromov-Witten invariants of the product stack $\mathcal{X}_1\times \mathcal{X}_2$ in terms of Gromov-Witten invariants of the factors $\mathcal{X}_1$ and $\mathcal{X}_2$. As an application, we deduce a decomposition result for Gromov-Witten theory of trivial gerbes.
\end{abstract}

\tableofcontents

\section{Introduction}
Orbifold Gromov-Witten theory was constructed by Chen and Ruan \cite{CRogw} for symplectic orbifolds and by Abramovich, Graber and Vistoli \cite{AGVaoqc}, \cite{AGV06} for smooth Deligne-Mumford stacks. Let $\mathcal{X}_1$ and $\mathcal{X}_2$ be two smooth proper Deligne-Mumford stacks with projective coarse moduli spaces. It is natural to ask if there is any relationship between the Gromov-Witten theory of $\mathcal{X}_1\times \mathcal{X}_2$ and the Gromov-Witten theories of $\mathcal{X}_1$ and $\mathcal{X}_2$. When the targets are smooth projective varieties, this question was answered by K. Behrend \cite{Beh97}, in which a product formula expressing Gromov-Witten invariants of $X_1\times X_2$ in terms of those of $X_1$ and $X_2$ is proven. 

In this paper we prove a formula, {\bf Theorem \ref{GW_class_product_formula}},  which expresses orbifold Gromov-Witten invariants of $\mathcal{X}_1\times \mathcal{X}_2$ in terms of orbifold Gromov-Witten invariants of $\mathcal{X}_1$ and $\mathcal{X}_2$. This formula generalizes the product formula of Behrend \cite{Beh97}. {\bf Theorem \ref{GW_class_product_formula}} is obtained as a consequence of some (technical) results concerning relationships between virtual fundamental classes, see {\bf Theorems \ref{virtual_class_formula} and \ref{weighted_virtual_class_formula}}. Our approach to the product formula in  orbifold Gromov-Witten theory closely follows that of \cite{Beh97} and makes extensive use of log geometry, especially Olsson's work on log twisted curves \cite{OLogCurv}, in order to handle various new phenomena which occur in the presence of stack structures.

Let $G$ be a finite group and $BG$ the classifying stack of $G$. The product $\mathcal{X}\times BG$ is the trivial $G$-gerbe over $\mathcal{X}$. Our results together with the results of Jarvis-Kimura \cite{JK} allow us to compute Gromov-Witten invariants of this trivial gerbe and verify the Gromov-Witten theoretic decomposition conjecture \cite{HHPS} in this most basic case, see Section \ref{trivial_gerbe}.

The rest of this paper is organized as follows. Section \ref{preparation} contains various preparatory definitions and results. The whole Section \ref{virtual_class_identity} is devoted to prove the main result {\bf Theorem \ref{virtual_class_formula}} on virtual fundamental classes. In Section \ref{GW_theory} we present the product formula in orbifold Gromov-Witten theory. Section \ref{trivial_gerbe} is devoted to the application to trivial gerbes.

\subsection{Notations and Conventions}
In this paper we fix the complex numbers as ground field.
By an  {\em algebraic stack}  we mean an  algebraic stack over $\bbC$   in the sense of \cite{Art74}.  By a {\em Deligne-Mumford stack} we mean an algebraic stack over $\bbC$ in the sense of \cite{DM69}. 
We assume moreover all stacks (and schemes) are  quasi-separated, locally noetherian, locally of finite type.
Following \cite{KatoLog}, logarithmic structures  are considered on the \'etale site of schemes. For the  extension of the theory to stacks, which will be assumed, see \cite{OlssLog03}. Stalks of sheaves are as in {\em \cite{Mil80}, page 60}. 

We use the curly letters ($\clC$)  to denote twisted curves and  capital  letters ($C$) to denote their coarse moduli spaces.
We denote by greek letters with  $\ \tilde{}\ $ gerby dual graphs and by greek letters without $\ \tilde{}\ $ the underlying dual graphs. Unless otherwise explicitly stated, the order of the isotropy groups of marking gerbes of twisted  curves will be denoted by the letter $b$, while the order of the stabilizer groups of nodes will be denoted by the greek letter $\gamma$ (cfr. {\bf Definition \ref{gerby_graph_def}}).

Let $\clX$ be a proper smooth Deligne-Mumford stack with projective coarse moduli space $\pi: \clX\to X$. The {\em inertia stack of $\clX$} may be defined as the fiber product over the diagonal
$$I(\clX):=\clX\times_{\clX\times \clX} \clX.$$ 
Alternatively $I(\clX)$ can be seen as the stack of representable morphisms from constant cyclotomic gerbes to $\clX$,
$$I(\clX)=\cup_{r\in \bbN} \text{HomRep}(B\mu_r, \clX).$$
This description yields a natural decomposition
$$I(\clX)=\cup_{r\in \bbN} I_{\mu_r}(\clX), \quad I_{\mu_r}(\clX):= \text{HomRep}(B\mu_r, \clX).$$
Clearly the stack $I_{\mu_r}(\clX)$ is non-empty only for finitely many $r\in \bbN$. 

The group $\mu_r$ acts on every object of $I_{\mu_r}(\clX)$. Removing this $\mu_r$ by rigidification (\cite{AOV08})  yields another stack 
$$\bar{I}(\clX):=\cup_{r\in \bbN} \bar{I}_{\mu_r}(\clX),$$
which is called the {\em rigidified inertia stack of $\clX$}. There are natural maps $I_{\mu_r}(\clX)\to \bar{I}_{\mu_r}(\clX)$ and $I(\clX)\to \bar{I}(\clX)$. Alternatively $\bar{I}(\clX)$ may be described as the stack of representable maps from (not necessarily trivial) cyclotomic gerbes to $\clX$. 

There are two naturally defined locally constant functions on $I(\clX)$:
$$\text{ord}: I(\clX)\to \bbN, \quad \text{age}: I(\clX)\to \mathbb{Q},$$
where $\text{ord}$ sends an object $[B\mu_r\to \clX]$ to the order $r\in \bbN$, and $\text{age}$ is the ``age'' function. These two functions are also defined for $\bar{I}(\clX)$, and by abuse of notation we denote them also by $\text{ord}$ and $\text{age}$. 

More details about (rigidified) inertia stacks and their properties can be found in e.g. {\em \cite{AGV06}}. 

We denote by $H_2^+(X, \bbZ)$ the semi-group of effective curve classes of $X$. 

We refer the readers to \cite{AGVaoqc} and \cite{AGV06} for the constructions of orbifold Gromov-Witten theory. We will write $\clK_{g,n}(\clX, \beta)$ for the moduli stack of twisted stable maps to $\clX$ with specified discrete data $g,n,\beta$. When the discrete data is clear from the context or irrelevant for the discussion, we simply denote the moduli stack by $\clK(\clX)$.

\subsection*{Acknowledgment}
We thank D. Abramovich, K. Behrend,  T. Graber, A. Kresch, F. Nironi, M. Olsson, and A. Vistoli for useful discussions. Y. J. and H.-H. T. thank Mathematical Sciences Research Institute, where part of this paper is written, for hospitality and support. H.-H. T. is supported in part by NSF grant DMS-0757722.

\section{Preparatory results}\label{preparation}
\subsection{Prestable curves}
We start by recalling the notion of {\em prestable curve} and of {\em morphism of  prestable curves}.  

\begin{definition}[\cite{BehMan}, Definition 2.1]
\label{prest_curv_def}
A {\em prestable curve} over a scheme $T$ is a flat proper morphism $\pi:C\to T$of schemes such that the geometric fibers of $\pi$ are reduced, connected and 1-dimensional and have at most ordinary double points (nodes)  as singularities.
The {\em genus} of a prestable curve  $C\to T$ is the map $t\mapsto H^1(C_t,\clO_{C_t})$, which is a locally constant function $g:T\to \bbZ_{\geq 0}$.
\end{definition}
\begin{definition}[\cite{BehMan}, Definition 2.1]
\label{prest_curv_morph_def}
A {\em morphism of prestable curves} $p:C\to D$ over $T$ is a $T$-morphism of schemes such that for every geometric point $t\in T$ we have 
\begin{enumerate}
\item if $\eta$ is the generic point of an irreducible component of $D_t$, then the fiber of $p_t$ over $\eta$ is a finite $\eta$-scheme of degree at most one;
\item if $C'$ is the normalization of an irreducible component of $C_t$, then $p_t(C')$ is a single point only if $C'$ is rational. 
\end{enumerate}  
\end{definition}
Prestable curves admit infinitesimal automorphisms. Therefore their moduli stack is an Artin stack. The algebraic stack of $n$-pointed, genus $g$ prestable curves (see \cite{Beh97GW}), usually denoted by $\frM_{g,n}$ is a  smooth algebraic stack of dimension $3g-3+n$. It is not separated nor finite type.

Imposing {\em  stability conditions} (\cite{DM69}), which amounts to requiring that every irreducible rational component of the  curve  has at least three special points, and that  any irreducible genus 1  component has  at least one special point, gives objects with finite isomorphisms groups.  
Given a prestable curve, there is a construction called {\em stabilization} which produces a stable curve,
see \cite{KnudII} for details. A relative version of the construction also exists, which we describe in {\em Section \ref{part-stab-sect}}.
\subsection{Gerby Modular Graphs}
We recall the definition of {\em dual  graph}\footnote{In this paper we do not use the full  categorical treatment of  graphs  developed in \cite{BehMan} (cfr.  \cite{Behagw99}).} of a prestable curve given in \cite{BehMan}.
\begin{definition}
\label{modular-graph-def}
A {\em graph} $\tau$ is a quadruple $(F_{\tau},V_{\tau},j_{\tau},\partial_{\tau})$, where $F_{\tau}$
and $V_{\tau}$ are finite sets, $\partial_{\tau}: F_{\tau}\to V_{\tau}$ is a map and $j_{\tau}:
F_{\tau}\to F_{\tau}$ an involution. We call $F_{\tau}$ the set of {\em flags}, $V_{\tau}$ the set of {\em vertices}, $S_{\tau}=\{f\in F_{\tau}|j_{\tau}f=f\}$ the set of {\em tails} and $E_{\tau}=\{\{f_{1},f_{2}\}\subset F_{\tau}|f_{2}=j_{\tau}f_{1}, f_1\neq f_2\}$ the set of {\em edges} of $\tau$. For $v\in V_{\tau}$, let
$F_{\tau}(v)=\partial_{\tau}^{-1}(v)$ and $|v|=\# F_{\tau}(v)$, the {\em valence} of $v$. 
\end{definition}
Graphs are drawn  by representing vertices by dots, edges by curves connecting vertices, tails by half open curves, connected only at their closed end to a vertex. In such a way we obtain a topological space which is called the {\em geometric realization} of the graph and is denoted by $\abs{\tau}$. We will always assume that $|\tau|$ is connected.

\begin{definition}[\cite{BehMan}, Definition 1.5]\label{modular-graph-def-II}
A {\em modular graph} is a graph $\tau$ endowed with a map $g_\tau:V_\tau\to \bbZ_{\geq 0}$; $v\mapsto g(v)$. The number $g(v)$ is called the {\em genus} of $v$. 
\end{definition}

\begin{definition}
The {\em Euler characteristic} of a graph $\tau$ is defined as
\begin{equation*}
\chi(\tau)=\chi(\abs{\tau})-\sum_{v \in V_{\tau}}g(v),
\end{equation*}
where $\abs{\tau}$ is the geometric realization of the graph. 
\end{definition}
\begin{definition}[\cite{BehMan}, Definition 1.6]
Let $A$ be  a semigroup with indecomposable zero. An {\em $A$-structure} 
on $\tau$ is a map $\alpha:V_{\tau}\to A$. A pair $(\tau,\alpha)$ is called a 
{\em modular graph with $A$-structure or an $A$-graph}.  A {\em marked graph} is a pair $(A,\tau)$ where $A$ is a semigroup with indecomposable zero  and $\tau$ is an $A$-graph. 
\end{definition}
There is a notion of {\em stability} for graphs, namely
\begin{definition}[\cite{BehMan}, Definition 1.9]
A vertex $v$ of an $A$-graph $(\tau,\alpha)$ is stable if either $\alpha(v)\neq 0$ or $\alpha(v)=0$ and $2g(v)+\abs{v}\geq 3$. The $A$-graph $\tau$ is  called {\em stable} if all of its vertices are stable.
\end{definition}
Whenever an $A$-graph   is not stable there is a canonical procedure
to ``stabilize'' it, which we describe below. Before that we need to recall some notions of morphisms of graphs
introduced in  \cite{BehMan}.
\begin{definition}[\cite{BehMan}, Definition 1.7]
\label{comb_mor_def}
Let $(\sigma,\alpha)$  and $(\tau,\beta)$ be $A$-graphs. 
A {\em combinatorial morphism} $a:(\sigma,\alpha)\to (\tau,\beta)$ is a pair of maps $a_F: F_\sigma\to F_\tau$ and $a_V: V_\sigma\to V_\tau$ satisfying the following conditions:
\begin{enumerate}
\item the diagram
\begin{equation*}
\xymatrix{
F_\sigma\ar[r]^{\partial_\sigma}\ar[d]_{a_F} & V_\sigma\ar[d]^{a_V} \\
F_\tau\ar[r]_{\partial_\tau} & V_\tau
}
\end{equation*}
is commutative. 
\item For every $v\in V_\sigma$ let $w=a_V(v)$. The  induced 
  map $a_{V,v}: F_\sigma(v)\to F_\tau(w)$ is injective;
\item Let $f\in F_\sigma$ and $\overline{f}=j_\sigma(f)$. If $f\neq \overline{f}$, there exists $n\geq 1$ and $2n$ not necessarily distinct flags $f_1,...,f_n,\overline{f}_1,...,\overline{f}_n\in F_\tau$ such that
\begin{enumerate}
\item[(a)] $f_1=a_F(f)$, $\overline{f}_n=a_F(\overline{f})$;
\item[(b)] $j_\tau(f_i)=\overline{f}_i$;
\item[(c)] $\partial_\tau(\overline{f}_i)=\partial_\tau(f_{i+1})$ for $i=1,..,n-1$;
\item[(d)] for all $i=1,...,n-1$ we have $$\overline{f}_i\neq f_{i+1}\Rightarrow g(v_i)\neq 0\quad \mbox{and}\quad \beta(v_i)=0;$$ 
\end{enumerate}
\item for every $v\in V_\sigma$ we have $\alpha(v)=\beta(a_V(v))$ and $g(v)=g(a_V(v))$.
\end{enumerate} 

A {\em combinatorial morphism of marked graphs} $(B,\sigma,\alpha)\to (A,\tau,\beta)$ is a pair $(\xi,a)$, where $\xi:A\to B$ is a homomorphism of semigroups and $a:(\sigma,\alpha)\to (\tau,\xi\circ\beta)$ is a combinatorial morphism of $B$-graphs. 
\end{definition}

According to {\em  \cite{BehMan}, Proposition 1.13},\label{graph-stab-page}  given an $A$-graph $\tau$ which is not stable, there exists a stable $A$-graph $\tau^s$, together with a combinatorial morphism $a: \tau^s\to\tau$, such that every combinatorial morphism $\sigma\to \tau$, where $\sigma$ is a stable $A$-graph, factors uniquely through $\tau^s$.  In [Ibid.], {\em Section 5} the following characterization  of the stabilization morphism is provided. For every edge of $\tau^s$  $\{f,\overline{f}\}$,  there exists a unique sequence of distinct edges  $\{f_1,\overline{f}_1\}$,...,  $\{f_n,\overline{f}_n\}$ in $E_\tau$ such that
$a(f)=f_1$, $a(\overline{f})=\overline{f}_n$, and $\partial\overline{f}_i=\partial f_{i+1}$ for all $i=1,...,n-1$. All vertices $v_i=\partial\overline{f}_i=\partial f_{i+1}$ have valence two. We call $\{f_1, \overline{f}_1,..., f_n, \overline{f}_n\}$ the {\em long edge} associated to $\{f, \overline{f}\}$. The edges $\{f_i, \overline{f}_i\}$ are called the {\em factors} of the long edge. 

For every tail $f$ of $\tau^s$ there exists a unique sequence of flags $f_1,...,f_n$ of $\tau$ such that
\begin{enumerate}
\item if $n$ is odd, then $\{f_1,f_2\}$,...,$\{f_{n-2},f_{n-1}\}$ are edges of $\tau$, $a(f_1)=f_1$, $\partial f_{2i}=\partial f_{2i+1}$, for $i=1,...,1/2(n-1)$, and $f_n$ is a tail of $\tau$;
\item if $n$ is even, then $\{f_1,f_2\}$,...,$\{f_{n-1},f_n\}$ are edges of $\tau$, $a(f)=f_1$,  $\partial f_{2i}=\partial f_{2i+1}$, for $i=1,...,n/2-1$, and $\partial f_n$ has valence one.
\end{enumerate}
 All vertices $v_i=\partial f_{2i}=\partial f_{2i+1}$, for $i=1,...,[\frac{1}{2}(n-1)]$ have valence two. We call $\{f_1,...,f_n\}$ the {\em long tail} associated to $f$. If $n$ is odd, we call the edges   $\{f_{2i-1},f_{2i}\}$, for  $i=1,...,1/2(n-1)$, and the tail $f_n$ the {\em factors} of this long tail.  If $n$ is even, we call the edges  $\{f_{2i-1},f_{2i}\}$ for  $i=1,...,n/2$ {\em factors} of the long tail.

Let $\clX$ be a smooth proper Deligne-Mumford stack with projective coarse moduli space. Orbifold Gromov-Witten theory of $\clX$ is based on the notion of twisted stable map developped in \cite{AV02}, \cite{AGVaoqc} and \cite{AGV06}.
The source curves of such maps are {\em twisted curves}, which we describe in the following. 

\subsection{Twisted curves}
In this section we recall the definition and some useful results about twisted curves which explain why  the definition of gerby dual graphs is well posed.

\begin{definition}[\cite{AV02}, Definition 4.1.2]\label{tw-curve-def}
Let $S$ be  a scheme. A {\em (balanced) twisted curve}  over $S$ is a proper flat Deligne-Mumford stack  $\clC\to S$ whose fibers are pure one-dimensional and geometrically connected with at most nodal singularities. Let $\clC\to C$ be the coarse moduli space of $\clC$, and let $C^{sm}\subset C$ be the open subset where $C\to S$ is smooth. Then the inverse image $\clC\times_C C^{sm}\subset \clC$ is equal to the open substack  of $\clC$ where $\clC\to S$ is smooth. Moreover for any geometric point $s\to S$ the map $\clC_s\to C_s$ is an isomorphism over some dense open subset of $C_s$. As a consequence the coarse moduli space is also a nodal curve. For any geometric point mapping to a node $p\to C$, there exists an \'etale neighborhood $\Spec(A)\to C$ of $p$ and an \'etale morphism $$\Spec{A}\to \Spec{(\clO_s[x,y]/xy-t)}$$
for some $t\in \clO_S$ such that the pullback $\clC\times_C \Spec{A}$ is isomorphic to
$$\left[\Spec{(A[z,w]/zw=t',z^n=x, w^n=y)}/\Gamma\right]$$
for some $t'\in \clO_S$, where $\Gamma$ is a finite cyclic group of order $n$ such that if $\gamma\in \Gamma$ is a generator then $\gamma(z)=\zeta\cdot z$ and  $\gamma(z)=\zeta^{-1}\cdot z$. This kind of action is called balanced.
A twisted curve has {\em genus} $g$ is the genus of $C_s$ is $g$ for any geometric point $s\to S$. An {\em $n$-pointed twisted curve} is a twisted curve together with  a collection of disjoint closed substacks $\{\Sigma_i\}_{i=1}^n$ such that
\begin{enumerate}
\item each $\Sigma_i$ is contained in the smooth locus of $\clC\to S$;
\item the stacks $\Sigma_i$ are \'etale gerbes over $S$;
\item if $\clC_{gen}$ denotes the complement of the $\Sigma_i$ in the smooth locus of $\clC\to S$ then $\clC_{gen}$ is a scheme.
\end{enumerate}
\end{definition}

A smooth  $n$-pointed prestable  twisted curve $\clC$  over an arbitrary scheme $S$
is determined up to unique isomorphism by its coarse moduli space $C$ and by the 
order of the isotropy groups of its marking sections. 
This is proven in \cite{AGV06}, {\em Theorem 4.2.1}.
\begin{theo}
\label{root_smooth_tw_curve}
Let $b_1,..,b_n\geq 1$ be integers numbers.
Let $C$ be an $n$-pointed smooth prestable curve over $S$.
 Let $S_i\subset C$ denote the image of the sections of $C\to S$. Let $L_i:=\clO_C(-S_i)$ and let $\sigma_i$ be the canonical sections vanishing on $S_i$. 
 Then there  is, up to isomorphisms, a unique  twisted curve $\clC$ over $S$ with coarse moduli space $C$ such that the $i$-th marked gerbe is banded by $\mu_{b_i}$.  It is obtained as
\begin{equation*}
\clC\simeq \sqrt[b_1]{(C,L_1,\sigma_1)}\times_C...\times_C  \sqrt[b_n]{(C,L_n,\sigma_n)}
\end{equation*}
where $\sqrt[b_i]{(C,L_i,\sigma_i)}$ denote the construction known as ``root of a line bundle with section'' \cite{Cadm03}, \cite{AGV06}.
\end{theo}

For families of nodal twisted curves there is no such a simple description in terms
of the coarse moduli space and of the order of the isotropy groups of the singular points.  As we will describe in {\em Section \ref{log-tw-curves-sect}}, logarithmic geometry is needed to 
describe properly singular twisted curves.  

In case the nodes are gerbes over the base of the curve, it  is still possible to describe a nodal twisted curve in terms of smooth twisted curves by means of a construction analogous to the {\em clutching construction} of \cite{KnudII}, extended  to  stacks in \cite{AGV06}, {\em Section 5} and {\em Appendix A}. However in general the  singular locus of a twisted curve is not a gerbe over the whole  base.  If the base of the curve has a non reduced structure this happens
also if the coarse moduli space curve has a generic node. 
 If we restrict to  twisted curves over the spectrum of an algebraically closed field,
we can always apply the construction of  \cite{AGV06} we alluded to above.
In fact the singular curve is isomorphic to the   pushout along
 smooth substacks   of  smooth twisted curves.
Over $\Spec{\bbC}$, given a prestable curve $C$, and given a set of integer numbers specifying the order of the isotropy groups of the nodes, there is only one twisted curve $\clC\to C$.  Uniqueness can be seen locally. Indeed, the singular locus can be described \'etale locally  as a quotient stack as  follows. For simplicity assume that $\clC$ has a single node of order $r$, and assume that \'etale locally near this node the coarse moduli space $C$ is isomorphic to $\Spec{\bbC[x,y]/(xy)}$.  According to \cite{OLogCurv}, an \'etale local description for $\clC$  near this node is
$$
[\Spec{\bbC[z,w]/(zw)/\mu_r}],
$$
where $\mu_r$ acts via $z\mapsto \zeta_r z, w\mapsto \zeta_r^{-1} w$, $\zeta_r:=\exp(2\pi \sqrt{-1}/r)$.

\begin{rema}
As we mentioned above, if the base has a non reduced structure, there can be non isomorphic twisted curves over $C$, even when the order of the  stabilizer group of the singular locus is assigned.   Let  $C$ be  a   the trivial family of  nodal curves, \'etale locally isomorphic to   
 $$\Spec{\bbC[x,y]/xy}\times \Spec{\bbC[\epsilon]/\epsilon^r}.$$ 
There are two families of twisted curves over $C$ with node of order $r$. One is the trivial family, the other  can be described 
 \'etale locally  as  
$$
[\Spec{(\bbC[\epsilon]/\epsilon^r)[z,w]/(zw-\epsilon)/\mu_r}],
$$
where $z^r=x$, $y^r=y$,  $\epsilon^r=0$, and the action of the cyclic group
given by the usual balanced action.
This  stack has a closed substack corresponding to a  $\mu_r$-gerbe over the closed point, but it does not contain a $\mu_r$ gerbe over the base. 
\end{rema}

We generalize  \textbf{Definitions \ref{modular-graph-def}} and \textbf{\ref{modular-graph-def-II}} to the notion of  {\em gerby modular graph}, which allows us  to  suitably stratify  stacks parametrizing twisted curves and twisted stable maps.

\begin{notation}\label{labeling_set}
Let $\mho$ be a finite set and $\mathfrak{o}: \mho\to \bbN$, $\mathfrak{a}: \mho\to \bbN\cup\{0\}$ two set maps.
\end{notation}

The following example is important. 
\begin{example}\label{labeling_set_from_stack}
Let $\clX$ be a smooth proper Deligne-Mumford stack. A triple $(\mho, \mathfrak{o}, \mathfrak{a})$ is obtained as follows. Let $\mho=\mho(\clX)$ be the set of connected components of the rigidified cyclotomic inertia stack of $\clX$. Let $\mathfrak{o}$ be the set map  $\mho(\clX) \to \bbN$ taking $U\in\mho(\clX)$ to the integer $r$ such that  $U\subseteq \overline{I}_{\mu_r}(\clX)$. Let  $\mathfrak{a}: \mho(\clX)\to \bbN\cup \{0\}$ be the set map taking $U$ to its age. 
\end{example}

\begin{definition}\label{gerby_graph_def}
Let $(\mho, \mathfrak{o}, \mathfrak{a})$ be as in {\bf Notation \ref{labeling_set}}. A {\em gerby modular graph} $\ttau$ associated to $(\mho, \mathfrak{o}, \mathfrak{a})$ is the data of an underlying  modular graph $\tau$ with a map $\mathfrak{g}:F_\ttau\to \mho$ such that $\mathfrak{g}(f)= \mathfrak{g}(f')$ whenever  the flags $f$, $f'$ form an edge  $\{f,f'\}\in E_\ttau$. We define $\gamma:=\mathfrak{o}\circ\mathfrak{g}$. 

Let $A$ be a semigroup with indecomposable zero.  A {\em gerby $A$-graph} is a gerby modular graph $\ttau$  whose underlying modular graph $\tau$ is endowed with an $A$-structure.
\end{definition}
\begin{rema}
We abuse the notation slightly: for a gerby modular graph $\ttau$ we denote the data $F_\tau, E_\tau$, e.t.c. associated to its underlying graph $\tau$ also by $F_\ttau, E_\ttau$, e.t.c.
\end{rema}

\begin{definition}\label{gerby_X_graph}
Let $\pi:\clX\to X$ be smooth proper Deligne-Mumford stack with projective coarse moduli space $X$, and let $(\mho=\mho(\clX), \mathfrak{o}, \mathfrak{a})$ be the triple obtained from $\clX$ as in {\bf Example \ref{labeling_set_from_stack}}. A gerby modular graph $\ttau$ associated to the triple $(\mho, \mathfrak{o}, \mathfrak{a})$, with $\mathfrak{g}: F_\ttau\to \mho(\clX)$ and an $H_2^+(X,\bbZ)$-structure $\alpha: V_\ttau\to H_2^+(X,\bbZ)$, is called a {\em gerby $\clX$-graph}. 
\end{definition}

The idea behind the definition of gerby $\clX$-graphs is as follows. Given a representable morphism $f: \clC\to \clX$ from a twisted curve $\clC$ over $\text{Spec}\, \bbC$ to a smooth Deligne-Mumford stack $\clX$, we would like to encode the structure of the domain curve $\clC$ combinatorially. The structure of the underlying coarse curve $C$ is encoded in a modular graph $\tau$ as in \cite{BehMan}. The stack structures on $\clC$ are recorded as follows. If $p\simeq B\mu_r \in \clC$ is a stack point, then the restriction $f|_p: B\mu_r\to \clX$ is a $\bbC$-point of the rigidified inertia stack $\bar{I}(\clX)$. The order of the stabilizer group of the stack point $p$ can be recovered from the value at the point $f|_p$ of the locally constant function $\text{ord}: \bar{I}(\clX)\to \bbN$.

\begin{definition}
The {\em dimension} of a gerby $\clX$-graph $\ttau$  is
\begin{equation*}
dim(\ttau,\clX):=\chi(\ttau)(dim\ \clX-3)-\int_{\beta(\ttau)}\omega_\clX + \# S_\ttau -\# E_\ttau -\sum_{f\in S_\ttau}\mathfrak{a}(\mathfrak{g}(f)),
\end{equation*}
where $\beta(\ttau):=\sum_{v\in V_\ttau} \alpha (v)\in H_2^+(X, \bbZ)$. 
\end{definition}

\subsection{Prestable curves and stable maps associated to gerby graphs}
As we mentioned in the introduction, gerby graphs will allow us
to stratify moduli spaces of curves and of maps according to  the
topological type  and the stacky structure of the curves. 
As in \cite{Beh97} the idea is that Gromov-Witten classes, which are classes  in $H_*(\overline{M}_{g,n})$,  are determined by their restriction  to the  boundary strata in  $\overline{M}_{g,n}$, which are labelled by modular graphs $\tau$. This restriction can be computed  in terms of intersection theory over   moduli spaces  of maps  whose  source curves have generically  the topological type and the stack structure specified by a   dual gerby graph $\ttau$ with underlying 
modular  graph $\tau$. We will compare the virtual fundamental classes of those 
moduli stacks for maps to $\clX_1\times \clX_2$ and to $\clX_i$. 

Let $\ttau$ be a gerby modular graph. We define {\em $\ttau$-marked prestable curve}, which generalizes Definition 2.6 of \cite{BehMan}.
\begin{definition}
A  {\em $\ttau$-marked prestable curve} over $T$ is a pair $(\clC,\Sigma)$, where $\clC=(\clC_v)_{v\in V_\ttau}$ is a family of twisted curves over $T$,   $\Sigma=(\Sigma_i)_{i\in F_\ttau}$ is a collection of disjoint \'etale gerbes banded by cyclic groups  over $T$ which are closed substacks $\Sigma_i\hookrightarrow \clC_v$
such that for every geometric point $t\in T$, 
\begin{enumerate}
\item $\Sigma_i|_t$ is contained in the smooth locus of $\clC_{\partial(i),t}$ for any $i\in F_\ttau$;
\item $\Sigma_i|_t$ and $\Sigma_j|_t$ are disjoint for $i\neq j$;
\item the coarse curve $C_{v,t}$ has genus $g(C_{v,t})=g(v)$ for all $v\in V_\ttau$;
\item $Aut(\Sigma_i|_t)\simeq \mu_r$ where $r=\gamma(i)$ for $i\in F_\ttau$.
\item for any $\{i,\overline{i}\}\in E_\ttau$ a band-inverting isomorphism  $\xi_{i\overline{i}}:\Sigma_i\isomto \Sigma_{\overline{i}}$ over $T$ is assigned, $\clC_{\partial(i)}$ and $\clC_{\partial(\overline{i})}$ are glued along $\Sigma_i$ and $\Sigma_{\overline{i}}$ via $\xi_{i\overline{i}}$.
\end{enumerate}
\end{definition}

The gerby graph $\ttau$ specifies how to glue the curves corresponding to the vertices in $V_\ttau$ along gerbes indexed by edges $\{i,\overline{i}\}$ in $E_\ttau$ via the isomorhisms inverting the band $\xi_{i\overline{i}}$. The twisted curves can be glued by doing a pushout in the category of Deligne-Mumford stacks. The gluing procedure associates to the node obtained after gluing a cyclic group of automorphisms leaving fixed the coarse moduli space of the curve. This is shown in the following Lemma.

\begin{lemm}[see \cite{ACV}, Proposition 7.1.1]\label{ghost-autom-lem}
Let $\clC_1$ and $\clC_2$ be two smooth twisted curves and let $i_1: \Sigma \hookrightarrow \clC_1$ and $i_2: \Sigma \hookrightarrow \clC_2$ be two $\mu_r$-gerbes. Let $\clC$ be the pushout along $\Sigma$. Then $Aut(\clC)_C\simeq \mu_r$, where $Aut(\clC)_C$ is the group of automorphisms of $\clC$ leaving $C$ fixed.
\end{lemm}

For a gerby modular graph $\ttau$, let $\mathfrak{M}^{tw}(\ttau)$ denote the moduli stack of $\ttau$-marked prestable curves. 

Let $\clX$ be a smooth proper Deligne-Mumford stack with projective coarse moduli space $X$. Let $\ttau$ be a gerby $\clX$-graph. Denote by $\alpha: V_\tau\to H_2^+(X, \bbZ)$ the $H_2^+(X, \bbZ)$-structure on the underlying modular graph $\tau$. For the notion of twisted stable map we refer to {\em \cite{AV02}, Definition 4.3.1}. 
\begin{definition}
A $\ttau$-marked twisted stable map over $T$ is a morphism $f:(\clC, \Sigma)\to \clX$ such that
\begin{enumerate}
\item
$(\clC=(\clC_v)_{v\in V_\ttau}, \Sigma=(\Sigma_i)_{i\in F_\ttau})$ is a $\ttau$-marked prestable twisted curve over $T$;
\item
$f:(\clC, \Sigma)\to \clX$ is a twisted stable map;
\item
For $i\in F_\ttau$, the restriction $f|_{\Sigma_i}$ is an object over $T$ of the component of $\bar{I}_\mu(\clX)$ indexed by $\mathfrak{g}(i)\in \mho(\clX)$; 
\item
Let $(C=(C_v)_{v\in V_\tau}, \overline{\Sigma}=(\overline{\Sigma}_i)_{i\in F_\tau})$ be the coarse curve and $\bar{f}: C\to X$ the induced map, then the data $(C, \overline{\Sigma}, \bar{f})$ is a stable $(X, \tau, \alpha)$-map over $T$ in the sense of {\em \cite{BehMan}, Definition 3.2}.
\end{enumerate}
\end{definition}
In particular, this means that $\bar{f}_*[C_v]=\alpha(v)$ for $v\in V_\tau$. 

Let $$\clK(\clX, \ttau)\subset \clK_{1-\chi(\tau), \#S_\tau}(\clX,\beta(\ttau))$$ be the moduli stack of $\ttau$-marked twisted stable maps to $\clX$. 

\subsection{Relative coarse moduli space}\label{rel_coarse_space}
We discuss the construction of relative coarse moduli spaces.

Let  $f:\clX\to\clY$ be a morphism of algebraic stacks locally of finite presentation. Let $I(\clX/\clY)=Ker(I(\clX)\to f^{*}\clY)=\clX\times_{\clX\times_\clY \clX}\clX$  denote the relative inertia stack. In [AOV08] the authors prove 
that there exists a relative coarse moduli space, namely
\begin{theo}
There exists an algebraic stack $X$, and morphisms $\clX\stackrel{\pi}{\to}X\stackrel{\overline{f}}{\to}\clY$ such that $f=\overline{f}\circ \pi$, satisfying the following properties:
\begin{enumerate}
\item $\overline{f}:X\to \clY$ is representable;
\item 
if there is $\clX\stackrel{\pi'}{\to}X'\stackrel{\overline{f}'}{\to}\clY$ with $\overline{f}'$ representable, then there is a unique $h:X\to X'$ such that $\pi'=h\circ \pi$ and $\overline{f}=\overline{f}'\circ h$;
\item $\pi$ is proper and quasi finite;
\item $\clO_X\to \pi_*\clO_\clX$ is an isomorphism;
\item if $X'$ is a stack and $X'\to X$ is a representable flat morphism, then $\clX\times_X X'\to X'\to\clY$ satisfies $(1)-(4)$ for the morphism $\clX\times_X X'\to\clY$.
\end{enumerate}
\end{theo}
If $\clX$ and $\clY$ are tame stacks over some scheme $T$, the formation of the relative coarse moduli space commutes with base change on $T$. This is the case of twisted stable maps we are considering. 
For  twisted curves over the complex numbers  there is the following local description of the relative coarse moduli space.   \'Etale locally, a twisted curve
is isomorphic to  $\clX=[V/\Gamma]$, where   $V$ is  be the spectrum of an henselian local ring and  $\Gamma$ is  a finite group acting on $V$ and  leaving fixed a geometric point 
$s\to V$.  Let $\clX\to\clY$ be a morphism of Deligne-Mumford type.  Let $K\subset \Gamma$ the subgroup fixing the composite $s\to\clX\to\clY$, namely $K$ is the kernel of the map $\Gamma\to Aut_\clY(f(s))$. Let $U=V/K$ be the geometric quotient, and $Q=\Gamma/K$. Let $\overline{\clX}=[U/Q]$. Then there exist a natural projection $q:\clX\to\overline{\clX}$ and a  unique factorization $g:\overline{\clX}\to\clY$ such that 
$f=g\circ q  $  identifying $\overline{\clX}$ as the relative coarse moduli space of $f$ (\cite{AOV08}, {\em Proposition 3.4}).
 \subsection{Partial stabilization}\label{part-stab-sect}
Let $\pi:C\to T$ be  a prestable curve and  let   $f:C\to M$  be a morphism from a prestable curve over $T$ to a smooth projective variety.  There is a canonical way to associate to it a stable map. This is proved in \cite{BehMan}. Let $N$ be a very ample invertible sheaf on $M$. Let  $L=\omega_{C/T}(\sum x_i)\otimes f^*N^{\otimes 3}$, where $x_i$ are the marked sections. The relative stabilization is defined as
$$
C'=\text{Proj}(\oplus_{k\geq 0}\pi_* (L^k)).
$$
The absolute stabilization is recovered from this construction by considering the case where  $M$ is a point. The {\em partial  stabilization} is obtained in the following way. Given a prestable curve $C$, suppose that it has a dual graph with some unstable vertex. Suppose we want to contract the components corresponding to some of those vertices. It is possible to choose a set of sections  makingthe curve stable. Leaving out the sections which make stable  the components to be contracted and stabilizing gives the partial stabilization. 

Given a morphism $f:\clC\to\clM$ from a twisted curve to a  smooth Deligne Mumford stack, \cite{AV02}, {\em Section 9} provides a construction which associates to $f$ a twisted stable map. If $f$ is not representable, one reduces to the representable case by taking the relative coarse moduli space as described above. Then one considers the morphism $f:C\to M$ induced by passing to the coarse moduli spaces. Let $C'\to M$ be the stable map obtained by taking the relative stabilization. Let $L'$ be a sufficiently relatively ample (in a sense that we will specify) line bundle with respect to $C'\to M$. Let $\clL$ be its pullback to $\clC$. Define $\clA:=\oplus_{k\geq 0}\pi_*(\clL)^k$, where $\pi:\clC\to\clM\times T$. The relative stabilization in the stack case is defined as 
\begin{equation*}
\clC':= \text{Proj}_\clM\ \clA\to\clM. 
\end{equation*} 
One verifies that $\clC'$ is flat over $T$ and $\clC'\to\clM$ is a twisted stable map.

\subsubsection{Dual graphs and stabilization}
We recall that,  as explained on page \pageref{graph-stab-page}, in \cite{BehMan}, {\em Proposition 1.13} the authors prove that, given an $A$-graph, where $A$ is a semigroup with indecomposable zero, there exists a stable $A$-graph $\tau^s$ and a combinatorial morphism $\tau^s\to \tau$ with the following universal property: any combinatorial morphism $\sigma\to\tau$, where $\sigma$ is a stable $A$-graph, factors uniquely through $\tau^s$. 

We adapt the notion of combinatorial morphism for gerby graphs.
\begin{definition}
\label{gerby_comb_mor_def}
A {\em combinatorial morphism of marked  gerby graphs} $$(B,\mho_B, \sigma)\to (A,\mho_A, \tau)$$ is a combinatorial morphism of the underlying marked graphs
$(a,\xi):(B,\sigma)\to (A, \tau)$, together with a set map $\omega: \mho_A\to\mho_B$  such that for all $f\in F_\sigma$,  $\mathfrak{g}(f)=\omega\circ \mathfrak{g}(a_F(f))$. 
\end{definition}
\begin{rema}
In this paper we will be concerned with the case where $A=H_2^+(X_1\times X_2, \bbZ)$, $B= H_2^+(X_i, \bbZ)$,  and $\mho_A=I(\clX_1\times\clX_2)$, $\mho_B=I(\clX_i)$.  
\end{rema}

\begin{example}\label{induced_map_between_labeling_sets}
Let $\clX$ and $\clY$ be two smooth proper Deligne-Mumford stacks with projective coarse moduli spaces $X$ and $Y$. Let $(\mho(\clX), \mathfrak{o}_\clX, \mathfrak{a}_\clX)$ and $(\mho(\clY), \mathfrak{o}_\clY, \mathfrak{a}_\clY)$ be the triples obtained from $\clX$ and $\clY$ as in {\bf Example \ref{labeling_set_from_stack}}. Let $F: \clX\to \clY$ be a  morphism. We can define part of the data constituting a combinatorial morphism of marked gerby graphs from $F$. More precisely we have a semi-group homomorphism $\overline{F}_*: H_2^+(X, \bbZ)\to H_2^+(Y, \bbZ)$, where $\overline{F}: X\to Y$ is the induced map between coarse moduli spaces. 

We define a set map $F_*: \mho(\clX)\to \mho(\clY)$ as follows. The morphism $F$ induces a morphism between the respective  cyclotomic inertia stacks $I(\clX)\to I(\clY)$ which maps connected components to connected components.  Indeed, let $\clG\to \clX $ be an object of $I_{\mu_r}(\clX)$ (or $\overline{I}_{\mu_r}(\clX)$). The morphism $\clG\to \clY$ obtained by composing with $F$ factors through the relative coarse moduli space gerbe (see Section \ref{rel_coarse_space}), which is isomorphic to $\clG$ if $F$ is representable. This defines canonically an object of $I_{\mu_s}(\clY)$ (or $\overline{I}_{\mu_r}(\clY)$) for some $s\in \bbN$ such that $s|r$. The order of the cyclic group banding the gerbe obtained as  relative coarse moduli space  divides the order of the cyclic group banding $\clG$.  
\end{example}

\subsection{Logarithmic geometry}
We give here a short introduction to {\em logarithmic geometry}, which provides a characterization of twisted curves we will  heavily use to prove some technical lemmas in the following.
Logarithmic structures have been introduced by Fontaine-Illusie and further studied by Kato \cite{KatoLog}. Given a scheme $X$, a {\em pre-logarithmic structure}, often called {\em pre-log structure},  consists of a constructible sheaf of monoids $M$ endowed with a morphism of monoids $\alpha:M\to \clO_X$, where    the structure sheaf is considered as a monoid with the multiplicative structure. Given a monoid or a sheaf of monoids $M$, we denote by $M^*$ the submonoid or the subsheaf of invertible elements.
 When the natural morphism $\alpha^{-1}(\clO_X^*)\to M^*$ is an isomorphism a pre-log structure is called a {\em log structure}. In this case we will indicate by  $\lambda: \clO_X^*\hookrightarrow M$ the natural inclusion of the subsheaf $\clO^*_X\subset \clO_X$ in the log structure $M$ over $X$. The quotient $M/\alpha^{-1}(\clO^*_X)$ is usually denoted by $\overline{M}$, and  called the {\em characteristic} or the {\em ghost sheaf}.   There is a canonical way to associate a log structure to a pre-log structure. Given a pre-log structure  $\alpha:M\to \clO_X$, the associated log structure, denoted $M^a$, is defined as the pushout in the category of sheaves of monoids as in the following diagram
\begin{equation*}
\xymatrix{
\alpha^{-1}(\clO^*_X)\ar[r]\ar[d]_{\alpha} & M\ar[d]\\
\clO^*_X\ar[r] & M^a.
}
\end{equation*}
The morphism to the sucture sheaf $\alpha^a:M^a\to\clO_X$ is induced by the pair of morphisms $(\alpha,\iota)$, where  $\iota:\clO^*_X\hookrightarrow\clO_X$ is the canonical inclusion. A scheme endowed with  a log structure $(X,M_X)$ is called a {\em log scheme}. Log schemes form a category. A morphism between two log schemes $(f,f^\flat): (X,M_X)\to (Y,M_Y)$ is a pair consisting of a morphism of schemes $f:X\to Y$ and a morphism of sheaves of monoids $f^\flat:f^*M_Y\to M_X$ compatible with the morphisms to the  structure sheave. The pullback of a log structure is defined
as the log structure associated to the pre-log structure obtained 
by taking the inverse image.\\

An important class of log structures are  {\em fine} log structures. By
definition this means that \'etale locally there exists  a morphism of sheaves of  monoids $\beta: P\to M$, where $P$ is the constant sheaf corresponding to a finitely generated  {\em integral}  monoid, such  that 
$M$ is isomorphic to the pre-log structure defined by $\alpha\circ \beta: P\to\clO_X$. A monoid $P$ is called integral if the cancellation law holds, namely, if the equality between elements of $P$ $x+y=x+y'$ implies $y=y'$.
  In this paper we will be concerned with  {\em locally free} log structures, namely fine log strutures such that for any geometric point $x\to X$ the stalk $\overline{M}_x$ is isomorphic to a free monoid $\bbN^r$ for some integer $r$. This turns out to be equivalent to the existence \'etale locally of a chart
given by a free monoid.

 Many basic constructions and definitions used in this paper are based on the notion of {\em simple morphism} of log structures.
A morphism between  free monoids $\phi: P_1\to P_2$ is called {\em simple} if $P_1$ and $P_2$ have the same rank, and for every irreducible element of $p_1\in P_1$ there exists a unique element $p_2\in P_2$ and an integer $b$ such that
$b\cdot p_2=\phi(p_1)$. A morphism of locally free log structures is called simple if it induces simple morphisms on the stalks. 

The simplest examples of log spaces are points. 
According to {\em \cite{KatoLog}, Example 2.5 (2)},  isomorphism classes of fine log structures over algebraically closed fields are in bijection with isomorphism classes of sharp integral monoids. A monoid is called {\em sharp} if it has no invertible element other than the unit element. Any integral   log structure over $X=\Spec{k}$
is of the form 
\begin{equation}\label{log_str_over_point}
M=\clO_X^*\oplus P  \mapsto \clO_X, \quad  (z,p)  \mapsto \begin{cases} z \quad \mbox{if}\quad  p=1\\
 0 \quad \mbox{if} \quad  p\neq 1.
 \end{cases}
\end{equation} 

\subsubsection{Automorphisms of logarithmic structures}\label{Automorphisms_Log_Str}
Throughout the paper we will often be concerned with automorphisms of locally free log structures. Here we make some remarks which will be useful in the following. Let $M$ be a locally free log structure over $X$, where $X$ is the spectrum of a henselian local ring. Denote by $x_0$ the closed point of $X$. We assume that $M$ admits a global chart $\bbN^r\to M$, inducing a bijection $\bbN^r\to \overline{M}_{x_0}$. Denote by $e_1,...,e_r$ the global sections of $M$ defined by the chart. Let $\phi:M\isomto M$ be an automorphism inducing the identity on $\overline{M}$.  Then $\phi$ corresponds to a collection of units $\vec{\zeta}=(\zeta_1,...,\zeta_{r})$ such that $\alpha(e_i)\zeta_i=\alpha(e_i)$. Let us first associate $\vec{\zeta}$ to $\phi$. For any $i=1,...,{r}$, $\phi(e_i)=e_i+\zeta_i$,  with  $\zeta_i\in \clO_X^*$. Indeed $\phi$ induces  the identity on $\overline{M}$, therefore $\phi(e_i)$ and $e_i$ can only differ by an unit. Since $M$ has integral stalks, the units $\vec{\zeta}$ are uniquely determined. On the other hand, a collection of units $\vec{\zeta}$ such that for any $i=1,..., r$, $\alpha(e_i)=\zeta_i\alpha(e_i)$ determines  an isomorphism of $M$. This is induced via the universal property of the pushout by the map $\bbN^r\to M$ mapping the $i$-th standard generator to $e_i+\zeta_i$ for $i=1,...,r$. Such an isomorphism  induces the identity on $\overline{M}$ and is  compatible with $\alpha$. 

We will use this equivalence between automorphisms and collection of units for locally free log structures over spectra of henselian rings, like discrete valutation rings, Artinian rings and  algebraically closed fields. 
\subsubsection{Log twisted curves}\label{log-tw-curves-sect}
From \textbf{Theorem \ref{root_smooth_tw_curve}} it follows that the stack
structure of the marking gerbes of a twisted curve is determined by taking the root of a   line bundle with section. As discussed thereafter, the stack structure at the nodes  cannot be described in this way, unless the nodes are gerbes over the base of the curve. In general, the natural language to describe  twisted nodes is log geometry \cite{OLogCurv}. 

 It is possible to show that given a family of prestable curves $C$ over $S$, there is a pair of log structures  $M_C$ over  $C$ and  $N_S$ over $S$, uniques up to unique isomorphism, which make the morphism of log schemes $(C,M_C)\to (S, N_S)$ {\em special}, in the sense of {\em Definition 2.6} of  \cite{OlssULS03} (see also {\bf Definition \ref{special_morphism}} below). This is proven in  [ibid.], {\em Theorems 2.7} and {\em 3.18}. Let $\clC$ be a twisted curve over $C$. The stack structure at the special points of  $\clC$  is determined by taking a ``root'' of the log structure over the base of $C$ in a sense that we will sketch below. Let us start by recalling how the canonical log structures are defined. A prestable curve $C\to S$ is an example of {\em semistable variety} in the sense of the following 
\begin{definition}
A scheme $X$ over a separably closed field $k$ is a {\em semistable variety} if for each closed point $x$ of $X$ there exists an \'etale neighborhood $(U,x')$ of $x$, integers $r\geq l$, and an \'etale morphism 
\begin{equation*}
U\to \text{Spec}\,k[X_1,..,X_r]/(X_1\cdot\cdot\cdot X_l).
\end{equation*}
\end{definition}
We need other preliminary definitions to be able to say what a special mophism of log schemes is.
\begin{definition}
A log smooth morphism $f:(X, M_X)\to (S, M_S)$ is {\em essentially semistable} if for each geometric point $x$ of $X$, the monoids $(f^{-1}\overline{M}_S)_{x}$
and $\overline{M}_{X,x}$ are free monoids, and if for suitable isomorphisms 
 $(f^{-1}\overline{M}_S)_{x}\simeq\bbN^r$,  $\overline{M}_{X,x}\simeq \bbN^{r+s}$ the map
\begin{equation*}
(f^{-1}\overline{M}_S)_{x}\to \overline{M}_{X,x}
\end{equation*}
 is of the form
\begin{eqnarray}
\label{ess_ss_eq}
e_i\mapsto \left\{ \begin{array}{cc}
e_i & \mbox{if}\ i\neq r,\\
e_r+...+e_{r+s} & \mbox{if}\ i=r.
\end{array}\right.
\end{eqnarray}
\end{definition}
If a log morphism $(X, M_X)\to (S, M_S)$ is essentially semistable locally in the smooth 
topology it admits the following description
\begin{lemm}
\label{ess_ss_lem}
Let $f:(X, M_X) \to (S, M_S)$ be an essentially semistable morphism of log schemes.   Then \'etale locally on $X$ and $S$ there exists charts $\bbN^r\to M_S$, 
$\bbN^{r+s}\to M_X$ such that the map $\bbN^r\to \bbN^{r+s}$ given by (\ref{ess_ss_eq}) is a chart for $f$, and such that the map
\begin{eqnarray}
\label{ess_ss_str_sh}
\clO_S\otimes_{\bbZ[\bbN^r]}\bbZ[\bbN^{r+s}]\to \clO_X
\end{eqnarray}
is smooth.
\end{lemm}
Let $f:(X,M_X)\to (S,M_S)$  be an essentially semistable morphism. Let $x\in X$ be a singular point. Then from   (\ref{ess_ss_str_sh}) we deduce that  in an \'etale neighborhood of $x$
\begin{equation*}
k \otimes_{\bbZ[\bbN^r]}\bbZ[\bbN^{r+s}]\simeq k[x_r,...,x_{r+s}]/(x_r\cdot\cdot\cdot x_{r+s} )\to \clO_X
\end{equation*}
is smooth. Therefore the morphism
\begin{equation*}
\overline{M}_{S,s}\to \overline{M}_{X,x}
\end{equation*}
is of the form  $\bbN^{r'}\to \bbN^{r'+s}$ with the map  as described in (\ref{ess_ss_eq}) for some $r'\leq r$.  This implies that there is only one irreducible element in the monoid $\overline{M}_{S,s}$ whose image in $\overline{M}_{X,x}$ is not irreducible. This defines a canonical map
\begin{equation*}
s_X:\{\mbox{singular points of}\  X\} \to \mbox{Irr}(\overline{M}_{S,s})
\end{equation*}
where we denote by Irr$(\overline{M}_{S,s})$ the set of irreducible elements of 
$ \overline{M}_{S,s}$. 
\begin{definition}\label{special_morphism}
\label{ess_ss_def}
An essentially semistable morphism of log schemes $f:X\to S$ is {\em special}
at a geometic point $s$ of $S$ if the map described above
\begin{equation*}
s_{X_s}:\{\mbox{singular points of}\  X_s\} \to \mbox{Irr}(\overline{M}_{S,s})
\end{equation*}
induces a bijection between the set of connected components of the singular locus of $X_s$ and Irr$(\overline{M}_{S,s})$.
\end{definition}
In  \cite{OlssULS03}, {\em Theorem 2.7} it is proven that if $(f, f'^\flat):(X, M'_X) \to (S,M'_S)$ is a  smooth, proper, integral\footnote{See {\em \cite{KatoLog}, Definition 4.1}. This implies flatness of the underlying morphism of schemes.}, vertical\footnote{This property means that  the cokernel of  $f'^\flat$ in the category of sheaves of monoids is a group.}  morphism of log schemes  and    $X$ is fiberwise  a semistable variety, then 
there exists a pair of log structures $(M_X,M_S)$ and a morphism $f^\flat: f^*M_S\to M_X$ making 
$(f, f^\flat): (X,M_X)\to (S,N_S)$ a special morphism. Moreover there are     morphisms of log structures 
\begin{equation*}
\phi:M_S\to M'_S,\quad\quad \psi:M_X\to M'_X
\end{equation*} 
and a cartesian diagram
\begin{equation*}
\xymatrix{
(X, M_X)\ar[d]_{(f,f^\flat)} \ar[r] & (X, M'_X)\ar[d]^{(f',f'^\flat)}\\
(S,M_S) \ar[r] & (S,M'_S).
}
\end{equation*}
The universal property of the log structures $(M_X, M_S)$ implies that,  given another pair $(\tilde{M}_X, \tilde{N}_S)$ of log structures making of $f:X\to S$ a special morphism, there are unique isomorphisms $\tilde{M}_X\isomto M_X$, 
$ \tilde{N}_S\isomto N_S$.

From the Theorem cited above it follows that the canonical  log structures associated to  a prestable curve or to a twisted curve are essentially uniques (up to unique isomorphism). Hence, to prove the existence, it is enough to define them \'etale locally.
Moreover, as observed in \cite{OlssULS03} in the proof of {\em Proposition 2.17},
due to  \cite{LMBca}, {\em Proposition 4.18}, there is an equivalence of categories
$$ \varinjlim_{U\to S} (\mbox{fine log structures on}\  X\times_S U)\to (\mbox{fine log structures on }\  X\times_S S_s)
$$
where $s$ is a point of $S$,   $S_s=\Spec{\clO^{sh}_{S,\overline{s}}}$ is  the spectrum of  the  (strict)  henselization of the local ring of $s$, $U\to S$ are \'etale neighborhoods of $s$. 
Then,  by a standard limit argument, 
 it is enough to   restrict  to curves over strictly henselian local rings.
Indeed, there exists some \'etale neighborhood $U$ such that
a given log structure over $S_s$ is obtained via pullback from  the corresponding log structure over $U$. 
 In this paper we work with stacks defined over fields of characteristics zero, therefore usually we make no distinction between strictly henselian and henselian rings. 

In the following we will denote by  $C$, resp. $\clC$, a prestable curve, resp. a twisted curve, over the spectrum of a henselian ring $S$. The canonical log structures associated to curves as above are the amalgamated sum (i.e. the pushout along $\clO_S^*$) of two log structures determined by  the nodes and by the marked points. Let us start by describing the log structure associated to the nodes. We need some further definitions.
Let $f:X\to S$ be a morphism of schemes. 
Let $t$ be a global section of $\clO_S$. Let $M_S$ be the log structure on $S$ associated to the pre-log structure $\bbN\to \clO_S$ mapping $1$ to $t$.
\begin{definition}[\cite{OlssULS03}]
A log smooth morphism $f:(X,M_X)\to (S,M_S)$ is {\em semi-stable} if for every geometric point $x\to X$ the stalk $\overline{M}_{X,x}$ is a free monoid and the map
\begin{equation*}
\bbN\to \overline{M}_{S,f(x)}\to\overline{M}_{X,x}
\end{equation*}
is the diagonal map. 
\end{definition}
The obstruction to the existence of a semistable log structure can be computed
by means of {\em Theorem 3.18} of  \cite{OlssULS03}. The obstruction space is given by the first \'etale cohomology group of a suitable subsheaf of $\clO_X^*$.
The canonical pair of log structures defined by  the singular loci of $C$ and $\clC$ are in fact semi-stable relative to $S$.  

Let us come to the definition for $f: \clC\to S$, following \cite{OLogCurv}. We will denote the  coarse moduli space  curve by $C$. 
Let $e_1,..,e_n$ be the nodes in the closed fiber of $C$. For $i=1,..,n$ let us choose open subsets $U_i$ containing only the node $e_i$. Let $\mathcal{U}_i$ be the preimage in $\clC$. Let $t_i\in\clO_S$ be an element such that \'etale locally $\clU_i$ is isomorphic to $$[\Spec{\clO_S[z,w]/(zw-t_i))}/\mu_{\gamma_i}],$$ where $\gamma_i$ is the order of the preimage of $e_i$ in $\clU_i$. Let $N_S^i$ be the log structure associated to the  pre-log structure $\bbN\to \clO_S$ mapping $1$ to $t_i$. It turns out that there exists a semi-stable log structure $M^i_{\clU_i}$ over $\clU_i$ relative to $(S,N_S^i)$.  This follows from an explicit computation of the obstruction group. This group vanishes only if the nodes are ``balanced''. Once a semi-stable log structure is given  over each $\clU_i$, then we have a global semi-stable log structure over $\clC$. 
Indeed, since $f^*N_S^i|_{\clU_i}\to M_{\clU_i}^i$ is semi-stable, it is an isomorphism out of the singular locus of $\clU_i$. If we denote by $\clU_i^-$ the complenet if $e_i$ in $\clC$, $\{\clU_i,\clU_i^-\}$ is an \'etale cover of $\clC$.
Over the intersection of the elements of the cover  $M_{\clU_i}^i$ and $f^*N_i$ glue  to a global log structure  $M^i_\clC$. Indeed  the restrictions are isomorphic and semistable log structures are unobstructed. 
 We define $N_S:=\oplus_{\clO_S^*} N_S^i$ and $M_\clC:=\oplus_{\clO_\clC^*} M_\clC^i$. The morphism $\clC\to S$ extends  to a semi-stable morphism 
\begin{equation*}
(\clC, M_\clC) \to (S, N_S).
\end{equation*}
As anticipated, there is also another log structure that we have to construct, which depends on the marked points of the curve. Let $\Sigma_i$ be the $i$-th marking gerbe. It corresponds to an invertible sheaf $\clI_i\to\clO_\clC$. From \cite{KatoLog}, {\em Complement 1} it follows that it defines a log structure $\clN_i$ over $\clC$. \'Etale locally it is possible to choose a generator  $f\in \clI_i$ for the ideal, and to define $\clN_i$ as the log structure associated to the pre-log structure $\bbN\to\clO_\clC$ mapping  $1$ to $f$. If $f'$ is another generator of $\clI_i$, corresponding to another log structure $\clN'_i$,  there exists a unique element $u\in \clO^*_\clC$ such that $uf'=f$. This unit defines an isomorphism $\clN_i\isomto \clN_i'$. Taking the amalgamated sum  $\tilde{M}_\clC:=M_\clC\oplus_{\clO_\clC^*} \clN_i$  (note that we use here a different notation with respect to  \cite{OLogCurv}) gives a global log structure such that the morphism $N_S\to  \tilde{M}_\clC$  is essentially semistable. 

The constructions defined above can be done of course for a prestable curve  with no stack structure. Let us consider in particular the coarse moduli space $C$ of $\clC$ and let us compare the canonical log structures  associated to $C$ with the log structures associated to $\clC$.
This can be made by looking  to the    \'etale local description. If the \'etale local description  around a node of the coarse moduli space is
$$
\Spec{A[x,y]/(xy-t)},
$$
where $A$ is an  henselian ring, then $\clC$  is isomorphic to
$$
[\Spec{A[z,w]/(zw-t')/\mu_r}],
$$
where $x=z^r$ and $y=w^r$, for some $t'\in A$ such that $t'^{r}=t$. 

Let us denote the  log structure over $S=\Spec{A}$ associated to $\clC$  by $N_S^\clC$ and the log structure associated to $C$ by $N_S^C$. We can see that $N_S^\clC$ can be thought of as a  root of $N_S^C$.
Indeed, $N_S^C$ is associated to the pre-log structure $\bbN\to A$ mapping
$1$ to $t$, while $N_S^\clC$ is associated to the pre-log structure $\bbN\to A$ mapping $1$ to $t'$. Because of that we see that  there is a morphism of pre-log structures
\begin{equation*}
\xymatrix{
\bbN\ar[rr]^{\times r }\ar[rd]& & \bbN\ar[ld]\\ 
& A & 
}
\end{equation*}
which induces a morphism of log structures 
\begin{eqnarray}
\label{log_str_clC_C}
N_S^{C}\to N_S^{\clC}.
\end{eqnarray}
Taking the characteristic and the stalk at the closed point $s_0$ of $\Spec{A}$ we get 
\begin{equation*}
\overline{N}_{S,s_0}^{C}\to \overline{N}_{S,s_0}^{\clC},
\end{equation*}
which is a simple  morphism of  monoids given again by
multiplication by $r$.   The log structure $N_S^{\clC}$ can be in fact  identified with 
the log structure $N'_S$ which is part of the data defining a log twisted curve 
over a prestable curve according to \textbf{Definition  \ref{log_tw_curve_def}}. 

 The comparison of the log structures associated to  the marked points of $C$  and  to the marking gerbes  of $\clC$ is analogous. However, here it is worth to notice that the log structures  $\clN_i$ associated  to the  marked points or to the marking gerbes  are equivalent to sets of line bundles with sections. This is proven in \cite{KatoLog}, where this kind of log structure are called {\em Deligne-Faltings} log structures. This is of course what one expects, since  by {\bf Theorem \ref{root_smooth_tw_curve}} the stack structure of the marking gerbes  is obtained by taking roots of line bundles  with sections over the coarse moduli space curve.

We now give the notion of  {\em n-pointed log twisted curve}  introduced in \cite{OLogCurv}, {\em Definition 1.7}
\begin{definition}\label{log_tw_curve_def}
An {\em $n$-pointed log twisted curve} over  a scheme $S$ is a  collection of data
\begin{equation*}
  (C/S,\{\sigma_i,a_i\}, l:N_S\to N'_S),
\end{equation*}
where $C/S$ is an n-pointed prestable curve, $\sigma_i:S\to C$ are sections, $a_i$, $i=1,...,n$ are integer-valued locally constant functions on $S$ such that for each $s\in S$ the integer $a_i(s)$ is positive and invertible over  $k(s)$, and $l:N_S\hookrightarrow N'_S$ is a  simple morphism of log structures over $S$, where  $N_S$ is the canonical log structure associated to $C/S$.
\end{definition}

Log twisted curves over $S$ form a groupoid. Let $(C^1/S,\{\sigma_i^1,a_i^1\}, l_1:N^1_S\hookrightarrow N'^1_S)$, $(C^2/S,\{\sigma_i^2,a_i^2\}, l_2:N^2_S\hookrightarrow N'^2_S)$ be two $n$-pointed log twisted curves. There are no morphisms between them unless $a^1_i=a^2_i$ $\forall i$. In this  case an isomoprhism is given by an isomorphism $\rho:  (C^1/S,\{\sigma_i^1\})\isomto  (C^2/S,\{\sigma_i^2\})$ of pointed prestable curves with an isomorphism $\epsilon:N'^1_S\isomto N'^2_S$ such that the diagram
\begin{eqnarray}
\label{log_curve_isom_diag}
\xymatrix{
N_S^1\ar[r]^{l_1}\ar[d]_{\simeq }^{\psi(\rho)} & N'^1_S\ar[d]^{\epsilon}\\ 
N^2_S\ar[r]^{l_2} & N'^2_S
}
\end{eqnarray}
commutes. Here  $\psi(\rho)$  is the isomorphism induced by the uniqueness of the canonical log structures   on $S$ and  $C$ (stated  in  \cite{OLogCurv},  {\em Theorem 3.6} and proved in \cite{OlssULS03},  {\em Proposition 2.7})  and by the isomorphism $\rho$, by means of which
$f^{*}_1 N^2_S\to \rho ^*M_{C_2}$ makes $f_1: C_1\to S$ a special morphism, whence the existence of a pair of isomorphisms $ N^1_S\isomto N_S^2$, $M^1_{C_1}\isomto \rho^* M^2_{C_2}$. 

\begin{rema}
\label{coarse_aut_noaffect_rem}
Let $f:C\to S$ be a prestable curve  and let $N_S$ be the associated canonical log structure over the base.   Let  $\rho:C\isomto C$ be an automorphism of prestable curves. By definition $\rho$ commutes with the projection $f$. Therefore 
 the induced automorphism $\rho^{*}N_S\isomto~N_S$ is the identity.
\end{rema}
\begin{rema}
\label{univ_prop_aut_iso_rem}
Let $N_S$ be a canonical log structure over the base of a prestable curve or of a twisted curve. Let us assume the curve has one node. There is a chart $\bbN\to N_S$, which defines a section that we denote by $e$.
 According to the discussion in {\em Section  \ref{Automorphisms_Log_Str}}, an automorphism of $N_S$ corresponds by the universal property of the pushout to a map $e\mapsto e+\lambda(u)$, where $\lambda:\clO_S^*\hookrightarrow N_S$ denotes the canonical inclusion. Let  $\epsilon: N^1_S\isomto N^2_S$ be a morphism between  two log structures as above induced by charts $\bbN\to S$ mapping $1$ to $t_1$, resp. $t_2$.   Let $e_1, e_2$ denote the sections of $N^1_S$, $N^2_S$  defined by the respective charts, then  $\epsilon(e_1)=e_2$. Again, since stalks are integral, and isomorphic to 
$\bbN\oplus \clO_{S,s_0}^*$, a unit $u$ such that $t_1 u=t_2$ is uniquely  associated to $\epsilon$. 
\end{rema}

\section{Identity for virtual fundamental classes}\label{virtual_class_identity}
The purpose of this Section is to prove {\bf Theorem \ref{virtual_class_formula}}, which is an identity relating the virtual fundamental class for moduli stack of twisted stable maps to $\clX_1\times \clX_2$ with virtual fundamental classes for moduli stacks of twisted stable maps to $\clX_1$ and to $\clX_2$.
\subsection{Set-up}\label{virtual_class_result_setup}
Let $\ttau$ be a gerby $(\clX_1\times \clX_2)$-graph with underlying modular graph $\tau$ and gerby structure $\mathfrak{g}: F_\ttau=F_\tau\to \mho(\clX_1\times \clX_2)$. The projections $p_i: \clX_1\times \clX_2\to \clX_i, i=1,2$ define maps $p_{i*}: H_2^+(X_1\times X_2, \bbZ)\to H_2^+(X_i, \bbZ)$. Let $\tau_i, i=1,2$ be the modular $H_2^+(X_i, \bbZ)$-graphs obtained as the stabilization of $\tau$ with respect to $p_{i*}$ (see {\em \cite{BehMan}, Remark 1.15}). By construction there are combinatorial morphisms $$(H_2^+(X_i, \bbZ), \tau_i)\to (H_2^+(X_1\times X_2, \bbZ), \tau), \quad i=1,2.$$ 
Consider the maps $F_{\tau_i}\to F_\tau$ between the sets of flags. The projections $p_i$ also define set maps $\mho(\clX_1\times \clX_2)\to \mho(\clX_i), i=1,2$ (see {\bf Example \ref{induced_map_between_labeling_sets}}). The compositions $$F_{\tau_i}\to F_\tau\overset{\mathfrak{g}}{\to}\mho(\clX_1\times \clX_2)\to \mho(\clX_i), \quad i=1,2,$$
define gerby structures on $\tau_i, i=1,2$. Denote by $\ttau_i, i=1,2$ the resulting $\clX_i$-graphs\footnote{It follows from the construction that if there exists a $\ttau$-marked twisted stable map to $\clX_1\times\clX_2$, then there exists $\ttau_i$-marked twisted stable maps to $\clX_i$, and there are combinatorial morphisms $\ttau_i\to \ttau$.}.
 
Let $\tau^s$ be the common absolute stabilization of $\tau, \tau_1, \tau_2$. 

Given an object $[f: (\clC/T, \Sigma)\to \clX_1\times \clX_2)]\in \clK(\clX_1\times \clX_2, \ttau)(T)$.  For $i=1,2$, applying relative coarse moduli space and relative stabilization constructions to the composite $\clC\to~ \clX_1\times \clX_2\times T\to \clX_i\times T$ yields representable morphisms $(\clC_i, \Sigma_i)\to \clX_i\times T$ fitting in the following diagram
\begin{equation}\label{inducing_maps}
\xymatrix{
(\clC, \Sigma) \ar[d]\ar[dr] &\\
(\clC_i, \Sigma_i) \ar[r] &\clX_i\times T.
}
\end{equation}
It is easy to check that $(\clC_i, \Sigma_i)\to \clX_i\times T$ in the above diagram is an object in $\clK(\clX_i, \ttau_i)$. This defines a morphism $$p=(p_1, p_2): \clK(\clX_1\times \clX_2, \ttau)\to \clK(\clX_1, \ttau_1)\times \clK(\clX_2, \ttau_2).$$

Throughout the rest this paper we will refer to the following diagram
\begin{eqnarray}\label{main_diag}
\xymatrix{
\clK(\clX_1\times\clX_2,\ttau)  \ar@/^/[rr]^{h}\ar@/^2pc/[rrr]^{p}\ar[r]\ar[d]& P'\ar[r]\ar@{}[rd]|{\square} \ar[d] & P\ar[r]_-{q}\ar[d]\ar@{}[rd]|{\square}& \clK(\clX_1,\ttau_1)\times \clK(\clX_2,\ttau_2)\ar[d]^{a}\\
 \frD^{tw}(\ttau) \ar@/_1pc/[rr]_l\ar[d]_{e}\ar[r]^{l_1} & \frD' \ar[r]^{l_2} &  \mathfrak{B}^{tw}\ar[r]^-{\phi}\ar[d]\ar@{}[rd]|{\square} & \mathfrak{M}^{tw}(\ttau_1)\times \mathfrak{M}^{tw}(\ttau_2)\ar[d]^{s\times s}\\
\mathfrak{M}^{tw}(\ttau) &  & \overline{M}(\tau^s)\ar[r]_-{\Delta}  &  \overline{M}(\tau^s)\times \overline{M}(\tau^s).
}
\end{eqnarray}
Some explanations are in order.
\begin{enumerate}
\item
$\overline{M}(\tau^s)$ is the moduli stack of $\tau^s$-marked stable curves. $\Delta$ is the diagonal morphism.
\item
As indicated in the diagram, $\frB^{tw}$ is the fiber product of $\overline{M}(\tau^s)$ and $\mathfrak{M}^{tw}(\ttau_1)\times \mathfrak{M}^{tw}(\ttau_2)$ over $\overline{M}(\tau^s)\times \overline{M}(\tau^s)$. 
\item
The morphism $\clK(\clX_i, \ttau_i)\to \frM^{tw}(\ttau_i)$ is obtained by forgetting the twisted stable maps but retaining the twisted curves.
\item
The stack $\frD'$ is defined to be $\frD(\tau) \times_\frB \frB^{tw}$ where $\frD(\tau)$ and $\frB$ are defined respectively in {\em \cite{Beh97}, diagrams (3) and (2)}. The map $l_2$ is the natural projection. 
\item 
The stacks $P$ and $P'$ are defined to be fiber products indicated in the diagram.
\end{enumerate}

 Diagram (\ref{main_diag}) is similar to {\em \cite{Beh97}, diagram (2)}. Other ingredients of (\ref{main_diag}) are explained in the rest of the Section.

\subsection{The stack $\frD^{tw}(\ttau)$}
We define $\frD^{tw}(\ttau)$ as the stack whose objects over $T$ are the following diagrams
\begin{eqnarray}\label{Dtwdiag}
 \xymatrix{
(\clC,\Sigma)\ar[r]^{\phi_1}\ar[d]^{\phi_2}& (\clC_1,\Sigma_1)\\
(\clC_2,\Sigma_2)
}
\end{eqnarray}
where
\begin{enumerate}
\item  $\clC$, resp. $\clC_i$,  is a $\ttau$-marked, resp. $\ttau_i$-marked,  twisted curve and  $\{\Sigma\}$, resp. $\{\Sigma_i\}$, denote a collection of disjoint \'etale gerbes in the smooth locus  over the base, of order specified by the gerby graph;  
\item  the morphisms $\phi_1$ and $\phi_2$  induce morphisms of prestable curves between the coarse moduli spaces;
\item the product morphism  $\phi_1\times\phi_2:\clC\to\clC_1\times\clC_2$ is {\em representable}, while $\phi_1$ and $\phi_2$ in general are not;   the morphisms $\phi_i:\clC\to\clC_i$ are {\em faithful};
\item the product morphism $\phi_1\times\phi_2$ does not admit infinitesimal automorphisms.  This implies that no unstable rational component of a fiber of $\clC$ is contracted in both $\clC_1$ and $\clC_2$;
\item  the orders of the isotropy group of any stack point of the curves
is  bounded from above by the largest order of the inertia groups of objects of $\clX_1\times \clX_2$. 
\end{enumerate}
\begin{rema}
The bound on the orders of the cyclic groups banding the marking gerbes ensures that the stack defined above is of finite  type over the stack $\frD'$.  This follows from   {\em \cite{OLogCurv}, Corollary 1.12} and  {\em \cite{Ols_HomSt}, Theorem 1.7} by applying an argument similar to {\em \cite{OLogCurv}, Theorem 1.16}. Moreover, again by {\em \cite{OLogCurv}, Corollary 1.12}, $\frD'$ is of finite type  
 over $\frD$ defined  in  \cite{Beh97}. As a consequence,   the stack relevant for our moduli problem  has a fundamental class (cfr. \cite{AGV06}). 
\end{rema}

\begin{rema}
Note that even though $\frD^{tw}(\ttau)$ parametrizes morphisms between stacks, it is in fact equivalent to a $1$-category.  Indeed  the morphisms  $\phi_1\times\phi_2: \clC\to\clC_1\times\clC_2$  are  representable and  map  a dense open representable substack  of $\clC$ to a dense  open representable substack of $\clC_1\times\clC_2$. According to e.g. {\em \cite{AV02}, Lemma 4.2.3}, any automorphism of such a morphism is trivial.
\end{rema}

The vertical arrow in diagram (\ref{inducing_maps}) defines the morphism $$\clK(\clX_1\times\clX_2, \ttau)\to \frD^{tw}(\ttau)$$ appearing in diagram (\ref{main_diag}).

There is a morphism $e: \frD^{tw}(\ttau)\to \frM^{tw}(\ttau)$ which takes an object (\ref{Dtwdiag}) to the $\ttau$-marked prestable twisted curve $(\clC, \Sigma)$. 
\begin{rema}
In what follows we will use the  valuative criterion for properness and the infinitesimal  lifting criterion for \'etaleness to prove properties of morphisms between stacks whose objects involve twisted curves. In those contexts we will not worry about the log structures associated to the marking gerbes. Indeed,  as far as those are concerned,  the criteria are satisfied because they are satisfied for the morphisms induced by passing to the coarse moduli spaces, which follows from \cite{Beh97}.
\end{rema}

\begin{prop}\label{e_is_etale}
  The morphism $e: \frD^{tw}(\ttau)\to \frM^{tw}(\ttau)$ is \'etale.
 \end{prop}
\begin{proo}
Since the stacks involved are algebraic stacks locally of finite type, it is enough to show that it is formally \'etale. We  use the infinitesimal lifting  criterion for \'etaleness (\cite{EGAIV_IV}, {\em Proposition 17.5.3 } and {\em Remark 17.5.4}, and {\em Definition 4.14} of \cite{LMBca}). We assume our stacks are locally n\"otherian, therefore when testing the criterion  we  can restrict to local  artinian  rings. 
Let $B$ be a local  artinian  ring, $I\subset B$ a square zero ideal and $A:=B/I$ the quotient ring. Consider an object in $\frD^{tw}(\ttau)(A)$:  
\begin{equation}\label{Dtw_object_over_A}
\xymatrix{
 \clC_A\ar[r]^{\phi_{1A}}\ar[d]^{\phi_{2_A}}  & \clC_{1A}\\
\clC_{2A}.
}
\end{equation}
We will show that (\ref{Dtw_object_over_A}) extends uniquely to $B$ once an extension $\clC_A\hookrightarrow \clC_B$ of $\clC_A$ to a twisted curve $\clC_B$ over $B$ is given. 

According  to \cite{Beh97}, {\em Lemma 4} we know that the coarse curves $C_{1A}$ and $C_{2A}$ extend uniquely to curves $C_{1B}$ and $C_{2B}$ over $B$. Recall that we have combinatorial morphisms $\tau_i\to\tau$, which include injective maps $a_V:V_{\tau_i}\to V_\tau$ for $i=1,2$. In order to show that $\phi_{iA}:C_A \to C_{iA}$ extend  to $B$ we  choose two sets of sections $S_i:= \{s^k_{iA}\}, i=1,2$ such that  for any $v\in V_\tau$ which is in the image of $a_V:V_{\tau_i}\to V_\tau$,  $S_i$ stabilizes the irreducible   components of  the curve $C_{v,0}$ not contracted in $C_{v_i,0}$, where $C_{v,0}$, $C_{v_i,0}$ denote  the  closed fibers.  This is possible because for any $\tau$-marked prestable curve $C$ over $\Spec{\bbC}$, there exists a graph $\tau'$ obtained from $\tau$ by adding some tails, such that $C$  is in the image of some point $C'$ of $\overline{M}(\tau')$ along the morphism
$\overline{M}(\tau')\to \frM(\tau)$. Moreover this morphism is smooth, therefore, given an extension of $C$ to a prestable curve $C_A$ over  some local artinian ring $A$, also  $C'$ extends to a stable curve  $C'_A$ such that $C_A$ is obtained from $C'_A$ by forgetting the sections corresponding to the additional tails of $\tau'$. By smoothness of $\overline{M}(\tau')\to \frM(\tau)$, one can extend the set of sections $S_i$ to $B$. To obtain the extensions $\phi_{iB}: C_B\to C_{iB}$ it is enough to forget the set of sections $S_j$, $j\neq i$, to stabilize and then forget the set of sections obtained as image of the sections $S_i$.

To show that twisted curves also extend, we  use the equivalence between twisted curves and log twisted curves  ({\bf Definition \ref{log_tw_curve_def}}). We observe that  the morphisms $\clC_A\to \clC_{iA}$ correspond to morphisms of the coarse moduli spaces together with  morphisms  of  the canonical  locally free log structures over the base  $A$ fitting in the following diagram,
\begin{equation*}
\xymatrix{
N_{iA}\ar@{^{(}->}[r]\ar[d]_{l_i} & N_A\ar[d]^{l}\ar[r]^{\alpha} & A\\  
N'_{iA}\ar@{^{(}->}[r] & N'_A\ar[r]^{\alpha'} & A.
}
\end{equation*}
We denote by $e_m$ the sections of $N$, $N'$ induced by the $m$-th standard generator of $\bbN^r$ via the charts $\bbN^r\to N$, $\bbN^r\to N'$, or equivalently
the corresponding edges in $E_\tau$, $E_\ttau$. (We assume here that the curve over the closed point has generic singularity type in order to use in the rest of the proof the language of dual graphs. The proof is identical if the curve has non generic singularity type). With this convention in the diagram above   $\alpha(e_m)=t_m$ and $\alpha'(e_m)=t'_m$, where $t_m, t'_m\in A$ and $t'^{\gamma_m}_m=t_m$, with  $\gamma_m=\gamma(e_m)$ being the order of the $m$-th node of $\clC$. Since we assume that we have an extension of $\clC$ to $B$,  an extension of the morphism $l:N_A\to N'_A$ to $B$  is given. This means that there are  given liftings of $t_m$, $t'_m$ to $B$, that we denote by $\tilde{t}_m$, $\tilde{t}'_m$. We now show that also $l_i:N_{iA}\to N'_{iA}$ extend to $B$  together with the morphisms $\phi_{iA}^{\flat}: N'_{iA}\to N'_{A}$. 
 Consider the diagram  of morphisms of pre-log structures over $\Spec{A}$
\begin{equation*}
\xymatrix{
\bbN^{r_i}_A\ar[d]_{(\vec{\gamma}_i)} \ar[r] & \bbN^r_A\ar[d]^{(\vec{\gamma})}\\ 
\bbN^{r_i}_A \ar[r]  & \bbN^r_A,
}
\end{equation*}
where the vertical arrows are diagonal matrices with entries the order of the isotropy  groups of the nodes of the closed fiber  corresponding to the irreducible elements of the monoids. (By construction  there are bijections  $\bbN^r\simeq \overline{N}_0$ and  $\bbN^{r_i}\simeq \overline{N}'_{i,0}$, where  $\overline{N}_0$, $\overline{N}'_{i,0}$  are the stalks at  the closed point.  The correspondence with generic nodes follows from {\bf Definition \ref{ess_ss_def}}).  The coarse moduli space morphisms $C\to C_i$ determine in fact the extensions of $l_i$ to $B$ and compatible morphisms $\phi_{iB}^\flat: N'_{iB}\to  N'_B$ as follows.   Let $J_{i,k}\subset\bbN$ be the  set of indices such that for $j_m\in J_{i,k}$, $e_{j_m}\in E_\tau$  partakes  in the long edge associated to $e_k\in E_{\tau_i}$.  The log structures $N'_{iB}$ determining  the stack structure of the nodes of any extension $\clC_{iB}$  compatible with $\clC_B$ must be  a sub log structure of $N'_B$. Therefore they must be the log structures associated to the pre-log structures defined as follows
\begin{equation*}
\oplus_{e_k\in E_{\ttau_i}}\bbN  \rightarrow B, \quad e_k  \mapsto \prod_{m\in J_{i,k}} {t'}^{\frac{\gamma_m}{\gamma_{ik}}}_m,
\end{equation*}
where $\gamma_m:=\gamma(e_m)$, $e_m\in E_\tau=E_\ttau$,   and $\gamma_{ik}:=\gamma(e_k)$, $e_k\in E_{\tau_i}=E_{\ttau_i}$.
This is the unique locally free log structure over $B$  which satisfies the required compatibility conditions.
\end{proo}

\subsection{The morphism $l$}
There is a morphism $$\frD^{tw}(\ttau)\to \frM^{tw}(\ttau_1)\times \frM^{tw}(\ttau_2)$$ which takes an object (\ref{Dtwdiag}) to the pair $((\clC_1, \Sigma_1), (\clC_2, \Sigma_2))$. There is another morphism $$\frD^{tw}(\ttau)\to \overline{M}(\tau^s)$$ which takes an object (\ref{Dtwdiag}) to the stabilization of the coarse curve $(C, \overline{\Sigma})$ of $(\clC, \Sigma)$. 

Since $\tau, \tau_1, \tau_2$ have the same absolute stabilization $\tau^s$, the following  diagram is 2-commutative,
\begin{equation*}
 \xymatrix{
& \frD^{tw}(\ttau)\ar[r]\ar[d]\ar@{}[rd]|{}   & \frM^{tw}(\ttau_1)\times \frM^{tw}(\ttau_2)\ar[d]^{s\times s}\ar[d]\\
& \overline{M}(\tau^s)\ar[r]^-{\Delta} & \overline{M}(\tau^s)\times \overline{M}(\tau^s),
}
\end{equation*}
which implies that there is a morphism $$l:\frD^{tw}(\ttau)\to \frB^{tw},$$ as in (\ref{main_diag}).

It follows easily from definitions that the morphism $l: \frD^{tw}(\ttau)\to \frB^{tw}$ factors through $l_1:\frD^{tw}(\ttau)\to  \frD'$, with the notation of diagram (\ref{main_diag}), i.e. $$l=l_2\circ l_1.$$ Recall that the stack $\frD'$ is defined as the fiber product 
$$\frD'=\frD(\tau)\times_{\frB}\frB^{tw}.$$
Unpacking the definitions of $\frB$ and $\frD(\tau)$ in \cite{Beh97}, we find that 
$$\frD'=\frD(\tau)\times_{\frM(\tau_1)\times\frM(\tau_2)}\frM^{tw}(\ttau_1)\times\frM^{tw}(\ttau_2). $$

\begin{rema}
\label{l_coarse}
The morphism $l$ induces a morphism
\begin{equation*}
\overline{l}: \mathfrak{D}(\tau)\to  \frB
\end{equation*}
between the stacks parametrizing the coarse moduli  spaces curves. This is the morphism $l$ in {\em \cite{Beh97}, diagram (2)}. We recall the proof in \cite{Beh97} that $\overline{l}$ is proper. Let $S_1, S_2$ be finite sets and $S:=S_1\coprod S_2$. Let $\tau_1'$ be a stable  modular graph obtained by adding $S_1$ to $S_{\tau_1}$ and $\tau_2'$  a stable  modular graph obtained by adding $S_2$ to $S_{\tau_2}$. Let $\tau'$ be the graph obtained by adding $S$ to $S_\tau$ in the unique compatible way in order for $\tau_i\to \tau$ to give  $\tau_i'\to \tau'$ inducing the inclusions $S_i\subset S$.  
We then have a cartesian diagram 
\begin{equation*}
\xymatrix{
\overline{M}(\tau')\ar[r]^-{\delta}\ar[d]\ar@{}[rd]|{\square} & \overline{M}(\tau_1')
\times \overline{M}(\tau_2')\ar[d]^{\chi}\\
\mathfrak{D}(\tau)\ar[r]^-{\tilde{\Delta}} & \mathfrak{M}(\tau_1)\times \mathfrak{M}(\tau_2)
}
\end{equation*}
where $\overline{M}(\tau')$,  $\overline{M}(\tau_1')$ and  $\overline{M}(\tau_2')$ are stacks of stable curves and the top row is a local presentation of $\tilde{\Delta}$. Notice that partial stabilization is used to define the top row.  Properness of $\overline{l}$ follows from properness of the induced morphism $\delta$ of the presentation. 
\end{rema}

In the rest of this Subsection we establish some properties of the morphism $l_1$, which will be important in proving {\bf Theorem \ref{virtual_class_formula}}.

\subsubsection{Properness}
\begin{prop}\label{l_1-proper-prop}
The morphism $l_1$ is proper.
\end{prop}
\begin{proo}
We apply the valuative criterion for properness. We work with stacks which are locally noetherian. Moreover the morphism $l_1$ is  of finite type. By {\em \cite{LMBca}, Proposition 7.8 and Theorem 7.10 (iii)}  we can apply  the valuative criterion  restricting  to  complete DVR's with algebraically closed  residue field. 
Let $R$ be such a DVR and let  $K$ be its field of fractions. Put $V=\Spec{R}$ and $U=\Spec{K}$. 

We first prove that $l_1$ is separated. Given two objects $\eta_1$ and $\eta_2$  in $\frD^{tw}(\ttau)(V)$, an isomorphism $\beta:l_1(\eta_1)\isomto l_1(\eta_2)$ in $\frD'(V)$ and an isomorphism  $\alpha:\eta_{1,U}\isomto\eta_{2,U}$  in $\frD^{tw}(\ttau)(U)$, we want to show that there exists a unique isomorphism $\tilde{\alpha}:\eta_1\isomto \eta_2$ restricting to $\alpha$ and such that $l_{1*}\tilde{\alpha}=\beta$.
This is not hard to see.
We only have to worry about isomorphisms associated to the stacky structure of 
$\clC$. We can equivalently prove the statement for automorphisms. Since $R$ is integral, a twisted curve carries automorphisms leaving fixed the coarse moduli space only if a node is a gerbe over $V$. In other words the pre-log structure associated to it maps to 0. Then $\clC_U$ and $\clC_V$ have the same automorphisms group. By definition of $\frD^{tw}(\ttau)$, $\phi_1\times\phi_2$ induces
an isomorphism between automorphisms of a node of $\clC$ and a subgroup
of the automorphism group of its image in $\clC_1\times\clC_2$. Under our hypotheses, $\beta$ encodes an isomorphism in this subgroup, hence $\alpha$ is determined uniquely. 

We make now the exercise of translating this argument in the language of logarithmic geometry.  
Consider two objects of $\frD^{tw}$, given by diagrams
\begin{equation*}
\xymatrix{
\clC^j \ar[r]^{(\phi^j_1, \phi^{\flat, j}_1)}  \ar[d]_{(\phi^j_2, \phi^{\flat, j}_2)} & \clC_1^j\\
\clC_2^j& 
}
\end{equation*}
where $j=1,2$. Here twisted curves are viewed as log twisted curves, and in particular as log schemes. In the following we denote by $\phi^{\flat, j}_i$ the morphisms induced between the canonical log structures over the bases of the curves, while we will not mention explicitly the log structures
over  the curves itselves. 

An isomorphism  between two objects  of  $\frD^{tw}(\ttau)(U)$ is  a triple of compatible isomorphisms $\vec{\rho}:=(\rho, \rho_1, \rho_2)$  of prestable curves, where $\rho:C^1\isomto C^2$, $\rho_i:C^1_i\isomto C^2_i$, $i=1,2$,  plus three compatible  isomorphisms  over $\vec{\rho}$ of  the simple morphisms of log structures defining the  log  twisted curves  $\clC^j$, $\clC_1^j$, $\clC^j_2$, $j=1,2$. The isomorphisms $\vec{\rho}$ between coarse curves have to satisfy the  compatibility conditions given by the commutativity of the following diagrams
\begin{eqnarray}\label{comp_rho}
\xymatrix{
C^1\ar[r]^{\phi_i}\ar[d]_{\sim}^{\rho}\ar@{}[rd]|{\curvearrowright} & C_i^1\ar[d]^{\sim}_{\rho_i}\\
C^2\ar[r]_{\phi_i}  & C_i^2
}
\end{eqnarray}
for $i=1,2$.
The isomorphisms of simple morphisms of log structures $\vec{\psi}=(\psi,\psi_1,\psi_2)$  have to satisfy the condition imposed by the commutativity of the diagrams
\begin{equation*}
\xymatrix{
&  & N^1\ar@{-->}[r]^{l^1}\ar@{-->}[dd]^-{\psi(\rho)} & N'^1\ar@{-->}[dd]^{\psi}\\
N_{i}^1\ar[urr]\ar[dd]_{\psi_i(\rho)}\ar[r]_{l^1_i}   & N'^1_i\ar[urr]\ar[dd]^{\psi_i} & \\
& & N^2\ar@{-->}[r]^{l^2} & N'^2\\
N_{i}^2\ar[urr]\ar[r]_{l^2_i}&   N'^2_i\ar[urr] & 
 }
\end{equation*}
for $i=1,2$, where the isomorphism $\psi(\rho)$, $\psi_i(\rho)$ are as  in (\ref{log_curve_isom_diag}) (see page \pageref{log_curve_isom_diag}).

Let $\bbN^r\to N^j$ and $\bbN^{r_i}\to N_i^j$, $j=1,2$, be the canonical log structures associated to $C^j$, $C_i^j$, with charts given by the pre-log structures defining them. In the following we will denote by  $e_m$ (resp. $e_{ik}$)
edges in $E_\tau$ (resp. $E_{\tau_i}$). Edges of the dual graph $\tau$ ( resp. $\tau_i$) correspond to irreducible elements of $\bbN^r$ (resp. $\bbN^{r_i}$) which are mapped to zero. We will use the symbols defined above decorated  by an upperscript $j$  to denote the sections of $N^j$, $N_i^j$  induced
by those irreducible elements. (cfr. {\bf Definition \ref{ess_ss_def}}).
  From {\bf Remark \ref{univ_prop_aut_iso_rem}} it follows that  the automorphisms $\psi(\rho)$ and $\psi_i(\rho_{i})$ correspond  to  collections of units in  $K^{*r}$ and $K^{*r_i}$, that we denote by $\vec{\eta}$ and by $\vec{\eta}_{i}$,  $i=1,2$, such that $\psi(\rho)(e_m^1)=e_m^2+\lambda(\eta_m)$,   $\psi_i(\rho_{i})(e_{ik}^1)=e_{ik}^2+\lambda_i(\eta_{ik})$.   Moreover, $\vec{\eta}$, $\vec{\eta}_{i}$ satisfy the following compatibility condition
\begin{equation*}
\prod_{m\in J_{i,k}} \eta_m= \eta_{ik},
\end{equation*}
for all $k$ such that $e_{ik}\in E_{\tau_i}$. Here  $J_{i,k}$ is the set of indices such that $m\in J_{i,k}$ if   $e_m\in E_{\tau}$ belongs to the long edge associated to $e_{ik}\in E_{\tau_i}$. The isomorphisms $\psi$, $\psi_i$ correspond to collections $\vec{\zeta}$,  $\vec{\zeta}_i$ in $K^{*r}$ and $K^{*r_i}$ which are  roots of elements of $\vec{\eta}$, $\vec{\eta}_i$. More precisely $\zeta_m$ is a $\gamma(e_m)$-th root of $\eta_m$, $\zeta_{ik}$ is a  $\gamma(e_{ik})$-th root of $\eta_{ik}$. For all $k$ such that $e_{ik}\in E_{\tau_i}=E_{\ttau_i}$ they  also  satisfy  the conditions
\begin{equation}
\label{psi_isom_Dtw_eq}
\prod_{m\in J_{i,k}} \zeta_m^{\epsilon^{ik}_m}=\zeta_{i,k},
\end{equation}
where $\epsilon^{ik}_m=\gamma(e_m)/\gamma(e_{ik})$, $e_m\in E_\ttau$.

The extension of the isomorphisms $\vec{\rho}$  to $V$ is 
encoded in $\beta$,  meaning in particular that $\vec{\eta}$ is in fact in $R^{*r}$, hence also $\vec{\zeta}$ is because $R$ is integrally closed in $K$.
 The  $\vec{\zeta}_i$ are also  in $R^{*r_i}$ by hypothesis. 
The compatibility  conditions  (\ref{psi_isom_Dtw_eq})
are  satisfied   in $R^{*r}\times R^{*r_i}$ since they are in $K^{*r}\times K^{*r_i}$. This concludes the proof that the morphism $l_1$ is separated. 

Since  $l_1$ is separated and  of finite type, we can apply  the valuative criterion for properness in the following form. Given an object of $\frD^{tw}(\ttau)(U)$,
\begin{eqnarray}\label{diag-over-U}
\xymatrix{
\clC_U\ar[r]\ar[d] & \clC_{1U}\\
\clC_{2U} &
}
\end{eqnarray}
and  an extension to $V$ of its image under $l_1$
\begin{equation}\label{part_stab_ext}
\clC_{1U}\hookrightarrow \clC_{1V}, \quad \clC_{2U}\hookrightarrow \clC_{2V},
\end{equation}
we have to show that  diagram (\ref{diag-over-U}) extends uniquely  to $V$  up to possible finite base extension. 

The diagram induced by  (\ref{diag-over-U}) by passing to the coarse moduli spaces  extends to $V$. The twisted curves over $V$ in (\ref{part_stab_ext}) are the data of   the coarse curves $C_{iV}$ plus simple morphisms of locally free  log structures
$l_i: N_{iV}\to N'_{iV}$. We recall that by construction those are the log structures associated to the pre-log structures defined in the following way.
Let $C_{i,0}$ be the closed fiber over $\Spec{R}$. We  will denote by the same symbol  $e_m$ a section of $N$ defined by an irreducible element of $\bbN^r$ via the chart $\bbN^r\to N$ and  the corresponding  node in the closed fiber of $C$.
 \'Etale locally around every node $e_m$  of $C_{i,0}$ the curve is isomorphic to 
$\Spec{R}[x,y]/(xy-t^{k_{i}^m\gamma_m})$. Here $k_{i}^m\in \bbZ_{\geq 0}\cup \{\infty\}$, and $\gamma_m:=\gamma(e_m)\in \bbN$ is the order of the stabilizer group of the node $e_m$ in $\clC_{i0}$ and $t$ is the uniformizing parameter. If $k_{i}^m=\infty$ we put $ t^{k_{i}^m}=0$.
 For any node $e_m$ let   $N_{iV}^{m}$ be  the  log structure associated to the pre-log structure
\begin{equation}
\label{def_eq_exp_Ni}
\bbN  \to R, \quad e  \mapsto t^{k_{i}^m\gamma_m},
\end{equation}  
and analogously let $N'^{m}_{iV}$ be  the  the log structure associated to  
\begin{equation}
\label{def_eq_exp_Ni'}
\bbN  \to R, \quad e  \mapsto t^{k_{i}^m}.
\end{equation}  
The canonical  log structures  $N_{iV}$ and  $N'_{iV}$ are defined as the amalgamated sums  $\oplus_{R^*} N_{iV}^{m}$  and   $\oplus_{R^*} N'^{m}_{iV}$. The simple morphisms  $l_i$ are induced by the morphisms of pre-log structures  
\begin{equation*}
\xymatrix{
\bbN^{r_i}_V\ar[r]^{(\vec{\gamma_i})}& \bbN^{r_i}_V
}
\end{equation*}
where $(\vec{\gamma}_i)$ is the matrix $diag(\gamma(e_1),...,\gamma(e_{r_i}))$ with $e_j\in E_{\ttau_i}$, $i=1,...,r_i$. 

We now  show that an extension $\clC_V$, $\phi_{i,V}$, $i=1,2$,  exists uniquely.   We  assume that the curves are not generically nodal because in this case the claim is trivial. Indeed, if $\clC_U$ has a generic node, the pre-log structure associated to that node is of the form  $\bbN\to K$ mapping 
$e$ to $0$, which extends obviously to $R$. 

  For $i=1,2$ we have the morphisms of log structures
\begin{equation*}
\xymatrix{
N_{i,V}\ar@{^{(}->}[d]\ar[r] & N'_{i,V}\\
N_V. & 
}
\end{equation*}
 The vertical arrow corresponds to the morphism of pre-log structures 
\begin{equation*}
\bbN^{r_i} \to \bbN^r, \quad e_k  \mapsto \sum_{m\in J_{i,k}} e_m,
\end{equation*}
 where $J_{i,k}$ is as  above the set of indices such that $\forall m\in J_{i,k}$, $e_m\in E_\tau$ is in the long edge associated to $e_k\in E_{\tau_i}$. 
For any $k\in \bbN$ such that $e_k\in E_{\tau_i}$, and for $i=1,2$, we have a  constraint  given by
\begin{equation*}
\sum_{m\in J_{i,k}}  k^m= k^k_i\cdot \gamma_{i,k},
\end{equation*}
where $k^m$, $k^k_i$  $\in \bbN_{\geq 0}\cup \{\infty\}$  are as in (\ref{def_eq_exp_Ni}) and  (\ref{def_eq_exp_Ni'}), and $\gamma_{i,k}:=\gamma(e_k)$ for $e_k\in E_{\ttau_i}=E_{\ttau_i}$.    Analogously, let $\gamma_j:=\gamma(e_j)$ for $e_j\in E_\ttau$. For all $j=1,...,\abs{E_\ttau}$ put $d_j:=gcd (k^j,\gamma_j)$ and $q_j:=\gamma_j/d_j$. Let $Q=max_j\{q_j\}$.  Let us consider the  finite  extension   $K'$ of $K$ obtained  by adding  the $Q$-th root of the uniformizing parameter, that we denote by $t'$. Let $R'$ be the integral closure of $R$ in $K'$. 
Over $V'=\Spec{R'}$ we can define the pre-log structure
\begin{equation*}
\bbN^r_{V'}  \to R', \quad e_j  \mapsto t^{k_j/\gamma_j}.
\end{equation*}
Notice that the  $t^{k_j/\gamma_j}$ is in $R'$.
We denote by $\tilde{N}'_{V'}$ the associated log structure. Let  $p:V'\to V$ denote  the morphism induced by the extension.  We have  a simple morphism of pre-log structures $$\tilde{l}: p^*N_V\to \tilde{N}_{V'}$$ giving  $p^*C_V$ the structure of a log twisted curve with the correct order of the stabilizers of the nodes, in a compatible way with the morphisms $p^*C_V\to p^*C_{iV}$. Moreover morphisms of log structures  $p^*N'_{iV}\to \tilde{N}'_{V'}$ are  also defined. This is the unique structure of log twisted curve which repects the prescription given by the fixed  morphisms of gerby dual graphs. 
\end{proo}

\begin{rema}
In the above proof the degenerate case corresponding to 
an unstable genus 1 component without marked points in $\clC$, contracted in either $\clC_1$ or $\clC_2$ does not appear. In this case the valuative criterion
is satisfied because of properness of $\overline{l}$. 
\end{rema}

\subsubsection{Orientation and degree calculation}

\begin{lemm}
The morphism $$\phi\circ l:\frD^{tw}(\ttau)\to \frM^{tw}(\ttau_1)\times \frM^{tw}(\ttau_2)$$ defines a natural orientation $[\phi\circ l]$.
\end{lemm}
\begin{proo}
The morphism $\phi\circ l$ is of relative Deligne-Mumford  type between smooth 
algebraic stacks of pure dimension. It  follows as observed in {\em \cite{Kre99}, footnote on page 529}
that it is an l.c.i. morphism, hence defines a natural orientation. 
\end{proo}

The following is a  {\em Proposition} in \cite{Beh97}
\begin{prop}
The morphisms $\Delta$ is a  proper regular local immersions. The orientations $[\Delta]$ and $[\phi]$ defined respectively by $\Delta$ and $\phi$ satisfy 
\begin{equation*}
[\phi]=(s\times s)^*[\Delta].
\end{equation*}
\end{prop}

\begin{proo}
 The morphism  $(s\times s)$ is flat, therefore $\phi$ is also a regular local immersion. The statement on orientations is obvious.
\end{proo}

\begin{prop}\label{prop:degree_of_l}
The morphism $l$ has degree $\mathfrak{c}$ where
\begin{eqnarray}\label{degree_of_l}
\mathfrak{c}:=\frac{\prod_{e_1\in E_{\ttau_1}} \gamma(e_1)\cdot \prod_{e_2\in E_{\ttau_2}} \gamma(e_2) }{\prod_{e\in E_\ttau} \gamma(e)}.
\end{eqnarray}
\end{prop}
\begin{proo}
We briefly  recall  the   definition of the stack of prestable  $\ttau$-marked twisted curves. The objects are collections of  prestable   twisted curves $\clC_v$ for any $v\in V_\ttau$, with $\abs{F_\ttau(v)}$ marked gerbes, such that the $i$-th gerbe corresponding to the flag $f_i\in F_\ttau(v)$  has order $\gamma(f_i)$. Moreover, as part of the defining data, there are given $\abs{E_\ttau}$ band-inverting isomorphisms of gerbes, prescribing how to glue gerbes corresponding to flags $\{f,\overline{f}\}\in E_\ttau$. 
While the order of the stabilizers of marking gerbes and of generic nodes, which are obtained by gluing marking gerbes corresponding to flags which partake in an edge, is fixed by the dual  gerby graph, this is not true for the order of the stabilizers of non-generic nodes. Moreover, since we are dealing with prestable curves, more nodes than those specified by the dual graph are allowed and those can have arbitrary cyclic groups as stabilizers. In order to get stacks of twisted curves of finite type over the corresponding stacks of prestable (coarse) curves, we constrain from above the order of the stabilizer of any special point. 

As before, we will consider  gerby  graphs $\ttau$ whose underlying modular graph admits a non-empty absolute stabilization. 

Generic curves are obtained by gluing smooth  twisted curves $\clC_v$, $v\in E_\ttau$ along marking gerbes as prescribed by the dual gerby graph $\ttau$. We can therefore assume that, given a curve over $\Spec{\bbC}$ corresponding to a generic point, no rational unstable component of $\clC$ is ever contracted in both $\clC_1$ and $\clC_2$. An automorphism of an object $(\vec{\clC}, \vec{\phi}):=(\clC,\clC_1,\clC_2,\phi_1,\phi_2)$ of $\frD^{tw}(\ttau)$ is given by  arrows and 2-arrows fitting  in the following 2-commutative diagram
\begin{equation*}
\xymatrix{
\clC\ar[r]^-{\phi_1\times \phi_2}\ar[d]_{\stackrel{\simeq}{\xi_\clC}} \ar@{}[rd]|{\stackrel{\Rightarrow}{\eta}}& \clC_1\times \clC_2 \ar[d]^{\stackrel{\simeq}{\xi_{\clC_1},{\xi_{\clC_2}}}} \\
\clC\ar[r]^-{\phi_1\times\phi_2} & \clC_1\times \clC_2. 
}
\end{equation*}
Notice however that the  2-arrows $\eta$ is  trivial in this case, because it is  an isomorphism between two representable  morphisms  of Deligne-Mumford stacks with trivial generic stabilizer (cfr. {\em \cite{AV02}, Lemma 4.2.3}).
Therefore automorphisms of $(\vec{\clC}, \vec{\phi})$ are given by compatible triples of automorphisms of twisted curves $(\xi, \xi_1,\xi_2)$. We will see in the following that for our  purposes it is convenient to consider the inertia 
of $\frD^{tw}(\ttau)$ relative to $\frD$. We only consider automorphisms leaving fixed the coarse moduli space. Since we required in the definition  {\bf (\ref{Dtwdiag})}  of $\frD^{tw}(\ttau)$ that the morphisms $\phi_i$ are faithful, we easily deduce that  automorphisms of $(\vec{\clC}, \vec{\phi})$ are determined by $\xi$.
For a twisted curve over $\Spec{\bbC}$ automorphism leaving the coarse moduli space fixed are as in {\bf Lemma \ref{ghost-autom-lem}}.
 As a consequence,  the generic stabilizer of $ \frD^{tw}(\ttau)$ is the following group
\begin{equation*}
\prod_{i=1}^r \mu_{\gamma_i},
\end{equation*}
where $\gamma_i$ is the order of the  $i$-th generic node and $r=\abs{E_\ttau}$.

The morphism $l_2$ in (\ref{main_diag}) is proper, since it is the base change of the morphism $\overline{l}$ which is proven to be proper in \cite{Beh97}. The morphism $l_1$ is proper as we proved in {\bf Proposition \ref{l_1-proper-prop}}. It is quasi-finite because e.g.  it fits in the commutative diagram
\begin{equation*}
\xymatrix{
\frD^{tw}(\ttau)\ar[r]^{l_1}\ar[d]^{e} & \frD'\ar[d]\\
\frM^{tw}(\ttau) \ar[r] & \frM(\tau),
}
\end{equation*}
where the left vertical arrow is \'etale and the  bottom horizontal arrow is quasi-finite (see. e.g. {\em \cite{AOV08}, Corollary A.8}). The stacks $\frD^{tw}(\ttau)$ and $\frD'$ are smooth stacks and $l_1$ induces
a bijection on geometric points. We have that the pushforward along $l_1$ of the orientation $[\frD^{tw}(\ttau)]$ is  a multiple by the degree of $l_1$  of the orientation  $[\frD']$.  The degree of a relative Deligne-Mumford type  morphism between Artin stacks can be computed locally in the   smooth topology of the base: Let $F:\clX\to\clY$ be  a generically quasi finite relative Deligne-Mumford morphism between Artin stacks, with $\clY$ integral.  Let   $U\to \clY$ be a  smooth surjective  morphism from a scheme, let $V$ the base change along $F$. The degree of $F$ is the degree of $V\to U$. 

Degrees of morphisms between  Deligne-Mumford stacks  can be computed  according to  \cite{Vist89}. 
Since the degree is multiplicative we can work relative to $\frD(\tau)$.  Both $\frD^{tw}(\ttau)$ and $\frD'$ are gerbes over the open (and dense) locus of $\frD(\tau)$ consisting of curves with only generic nodes. These gerbe structures correspond to the generic nodes of $\clC$ and respectively of $(\clC_1,\clC_2)$. This follows from the description of the automorphisms of objects of $\frD^{tw}(\ttau)$. Being a gerbe is preserved by base-change. 
The following cartesian  diagram  is 
obtained as a  base change of a smooth cover $D\to \frD(\tau)$ 
\begin{eqnarray}
\xymatrix{
D^{tw}\ar[r]\ar[d]\ar@{}[rd]|{\square} & D'\ar[r]\ar[d]\ar@{}[rd]|{\square} &  D\ar[d]\\
\frD^{tw}(\ttau) \ar[r] & \frD'\ar[r] & \frD(\tau).
}\nonumber
\end{eqnarray}
Both $D^{tw}$ and  $D'$ are generically Deligne-Mumford gerbes over $D$.
The degree of the structure morphisms of these generic gerbes can becomputed  from the relative automorphisms groups of objects of $\frD^{tw}(\ttau)\to \frD(\tau)$, $\frD^{tw}(\ttau)\to\frD'$, and $\frD^{tw}(\ttau) \to \frD(\tau)$.
We deduce 
\begin{equation*}
 deg(l_1)=\frac{\text{order}(\delta(\frD'\to \frD(\tau)))}{\text{order}(\delta(\frD^{tw}(\ttau)\to\frD(\tau)))},
 \end{equation*}
where $\delta$ denotes the relative inertia.
This proves the desired result, because $\frD'$ is birational to $\frB^{tw}$ (this follows from the corresponding statement for $\frD(\tau)$ and $\frB$ proven in \cite{Beh97}). 
\end{proo}

\subsection {Relative obstruction theories}
We compare relative obstruction theories on the stacks $\clK(\clX_1\times \clX_2, \ttau)$, $\clK(\clX_1, \ttau_1)$, and $\clK(\clX_2, \ttau_2)$. Consequently we prove {\bf Theorem \ref{virtual_class_formula}}.

\begin{prop}
 The following square is cartesian
\begin{equation}\label{sub_main_diag}
\xymatrix{
 \clK(\clX_1\times\clX_2,\ttau)\ar[r]^-{h}\ar[d]& P\ar[r]^-{q}\ar[d]\ar@{}[rd]|{\square}& \clK(\clX_1,\ttau_1)\times \clK(\clX_2,\ttau_2)\ar[d]^{a}\\
 \frD^{tw}(\ttau) \ar[r]^{l} & \mathfrak{B}^{tw}\ar[r]^-{\phi}& \mathfrak{M}^{tw}(\ttau_1)\times \mathfrak{M}^{tw}(\ttau_2).\\
}
\end{equation}
\end{prop}
\begin{proo}
The morphism $h:\clK(\clX_1\times\clX_2,\ttau)\to P$ is defined. Indeed we have a morphism $p_i: \clK(\clX_1\times\clX_2,\ttau)\to \clK(\clX_i,\ttau_i)$, $i=1,2$,  which is given by composing with the  projection to $\clX_i$ and  by  taking first the relative coarse moduli space and after the partial relative stabilization with respect to the morphism to $\clX_i$. The morphism $\clK(\clX_1\times\clX_2,\ttau)\to \frD^{tw}(\ttau)$ is given by constructing the induced morphisms to $\clK(\clX_i,\ttau_i)$   as described above and then forgetting all the maps to the target stacks and retaining
only the morphisms among the source curves. This gives a morphism $$\Theta: \clK(\clX_1\times\clX_2,\ttau)\to \frD^{tw}(\ttau)\times_{\frB^{tw}} P.$$

We now construct a morphism $\Xi: \frD^{tw}(\ttau)\times_{\frB^{tw}} P\to \clK(\clX_1\times\clX_2,\ttau)$. Given $f_1:\clD_1\to\clX_1$ , $f_2:\clD_2\to \clX_2$, an object $(\clC,\clC_1,\clC_2,\phi_1,\phi_2)$ of $\frD^{tw}(\ttau)$, and a pair of isomorphisms $\alpha_i:\clC_i\isomto \clD_i, i=1,2,$ we want to construct a representable morphism $f:\clC\to\clX_1\times\clX_2$. We claim that
\begin{equation}\label{reconstr_map_to_prod}
\xymatrix{
\clC\ar[r]^-{(\phi_1, \phi_2)} &  \clC_1\times\clC_2 \ar[rrr]^{(f_1\circ\alpha_1\circ pr_1,f_2\circ\alpha_2\circ pr_2)} & & &\clX_1\times\clX_2
}
\end{equation}
is a twisted stable map, namely  the composed  morphism in equation  (\ref{reconstr_map_to_prod}) is representable and stable..
This follows immediately from the definition of $\frD^{tw}(\ttau)$,  which parametrizes morphisms $[\phi_1:\clC\to\clC_1, \phi_2:\clC\to\clC_2]$ such that
$\phi_1\times \phi_2$ is representable and does  not admit infinitesimal
automorphisms relative to $\clC_1\times\clC_2$.

We want to prove that $\Theta\circ\Xi$ is isomorphic to the identity functor of the fiber product and that $\Xi\circ \Theta$ is isomorphic to the identity functor of $\clK(\clX_1\times\clX_2,\ttau)$.
It is clear that $\Xi\circ\Theta$ is the identity. It is also easy to check that $\Theta\circ\Xi$ is isomorphic to the identity. The key point is that due to the universal property of the relative stabilization we have morphisms $\clC_i^s\to \clC_i$, where $\clC_i^s$ is the curve obtained by applying   the relative coarse moduli space and the relative stabilization construction  to  the morphism  $\clC\to \clX_i$  given by the composition of  (\ref{reconstr_map_to_prod}) with the projection.  Since stable maps form a groupoid this morphism is in fact an isomorphism. This allows to define a natural transformation with the identity functor by composing arrows and 2-arrows arising 
canonically from the relative stabilization process. 
\end{proo}

Consider the universal family over $\clK(\clX_1\times \clX_2, \ttau)$:
\begin{equation*}
\xymatrix{
\clC\ar[r]^{f}\ar[d]^{\pi} & \clX_1\times \clX_2\\
\clK(\clX_1\times \clX_2, \ttau),
}
\end{equation*}
and the universal families over $\clK(\clX_i, \ttau_i), i=1,2$:
\begin{equation*}
\xymatrix{
\clC_i\ar[r]^{f_i}\ar[d]^{\pi_i} & \clX_i\\
\clK(\clX_i, \ttau_i).
}
\end{equation*}
According to \cite{AGV06} there are perfect relative obstruction theories
\begin{equation}\label{obs_theories}
\begin{split}
&E_{\clX_1\times \clX_2}^{\bullet\vee}:=R\pi_*(f^*\Omega_{\clX_1\times\clX_2}\otimes \omega_\pi)\to L_{\clK(\clX_1\times\clX_2, \ttau)/\frM^{tw}(\ttau)}\\
&E_{\clX_i}^{\bullet\vee}:=R\pi_{i*}(f_i^*\Omega_{\clX_i}\otimes \omega_{\pi_i})\to L_{\clK(\clX_i, \ttau_i)/\frM^{tw}(\ttau_i)}, \quad i=1,2.
\end{split}
\end{equation}
By {\bf Proposition \ref{e_is_etale}} the morphism $e: \frD^{tw}(\ttau)\to \frM^{tw}(\ttau)$ is \'etale. Hence $E_{\clX_1\times \clX_2}^{\bullet\vee}$ can be viewed as a relative perfect obstruction theory over $\frD^{tw}(\ttau)$. 

\begin{prop}\label{compare_obs_theories}
The two obstruction theories 
$l^*\phi^*(E_{\clX_1}^{\bullet}\boxplus E_{\clX_2}^{\bullet})$ and $E^{\bullet}_{\clX_1\times\clX_2}$ are naturally isomorphic.
\end{prop}
\begin{proo}
Consider the morphism 
$$p=(p_1, p_2): \clK(\clX_1\times \clX_2, \ttau)\to \clK(\clX_1, \ttau_1)\times \clK(\clX_2, \ttau_2).$$
Denote by $p_1^*\clC_1$ and $p_2^*\clC_2$ the pullback of the universal families over $\clK(\clX_1,\ttau_1)$ and $\clK(\clX_2,\ttau_2)$ to $\clK(\clX_1\times\clX_2,\ttau)$. This yields the following cartesian diagrams
\begin{equation*}
\xymatrix{
p_i^*\clC_i\ar[r]\ar[d]^{\pi_i'} \ar@/^/[rr]^{f_i'}\ar[d] & \clC_i\ar[r]_{f_i} \ar[d]^{\pi_i} & \clX_i\\
\clK(\clX_1\times\clX_2, \ttau)\ar[r]^{p_i} & \clK(\clX_i, \ttau_i).& 
}
\end{equation*}
Since this square is cartesian and the twisted curves are flat over their bases, we have (see e.g. \cite{BBrRuip09}, {\em Proposition A.85})
\begin{equation*}
p^* (R\pi_{1*}f_1^* T_{\clX_1}\oplus R\pi_{2*} f_{2}^*  T_{\clX_2})\simeq R\pi'_{1*} (f_1')^*T_{\clX_1}\oplus  R\pi'_{2*}(f_2')^*T_{\clX_2}.
\end{equation*}

Moreover we have over $\clK(\clX_1\times\clX_2, \ttau)$  morphisms $q_i: \clC\to p_i^*\clC_i$  which correspond to the relative stabilization morphism. Since for such morphisms we have $q_{i*}\clO_{\clC}\simeq \clO_{p_i^*\clC_i}$ we see that the functor $Rq_{i*}q_i^*$ is the identity.
Therefore 
\begin{equation}
\begin{split}
&R\pi'_{1*} (f_1')^*T_{\clX_1}\oplus  R\pi'_{2*}(f_2')^*T_{\clX_2}\\
\simeq\nonumber &R\pi'_{1*}Rq_{1*}q_1^*(f'_1)^* T_{\clX_1}\oplus R\pi'_{2*} Rq_{2*}q_2^*(f'_2)^*  T_{\clX_2}\\
\simeq\nonumber &R\pi_{*}f^*(T_{\clX_1}\boxplus T_{\clX_2}).
\end{split}
\end{equation}
The result follows.
\end{proo}

\begin{theo}\label{virtual_class_formula}
Let $[\clK(\clX_1\times\clX_2, \ttau)]^{vir}$ and $[\clK(\clX_i, \ttau_i)]^{vir}, i=1,2,$ be the virtual fundamental classes associated to the obstruction theories (\ref{obs_theories}). Then
\begin{equation}
h_*[\clK(\clX_1\times\clX_2, \ttau)]^{vir}=\mathfrak{c}\Delta^!([\clK(\clX_1, \ttau_1)]^{vir}\times [\clK(\clX_2, \ttau_2)]^{vir}).
\end{equation}
\end{theo}
\begin{proo}
Consider the diagram (\ref{main_diag}) and its portion (\ref{sub_main_diag}). We calculate 
\begin{equation*}
\begin{split}
&[\clK(\clX_1\times\clX_2, \ttau)]^{vir}\\
=&[\clK(\clX_1\times\clX_2,\ttau),E^\bullet_{\clX_1\times\clX_2}]\\
=& [\clK(\clX_1\times\clX_2,\ttau), l^*\phi^*(E^\bullet_{\clX_1}\boxplus E^\bullet_{\clX_2})] \quad (\text{by Proposition }\ref{compare_obs_theories})\\
=&(\phi\circ l)^![\clK(\clX_1,\ttau_1)\times  \clK(\clX_2,\ttau_2), E^\bullet_{\clX_1}\boxplus E^\bullet_{\clX_2}] \quad (\text{by \cite{BehFan}, Proposition 7.2})\\
=& (\phi\circ l)^!([\clK(\clX_1,\ttau_1),E^\bullet_{\clX_1}]\times[\clK(\clX_2,\ttau_2),E^\bullet_{\clX_2}]) \quad (\text{by \cite{BehFan}, Proposition 7.4})\\
=&
(\phi\circ l)^!([\clK(\clX_1, \ttau_1)]^{vir}\times [\clK(\clX_2, \ttau_2)]^{vir}).
\end{split}
\end{equation*}
One therefore has 
\begin{equation*}
\begin{split}
&\mathfrak{c}\cdot \Delta^!([\clK(\clX_1, \ttau_1)]^{vir}\times [\clK(\clX_2, \ttau_2)]^{vir})\\
=&\mathfrak{c} \cdot a^*(s\times s)^*[\Delta]([\clK(\clX_1, \ttau_1)]^{vir}\times [\clK(\clX_2, \ttau_2)]^{vir})\\
=&\mathfrak{c} \cdot a^*[\phi]([\clK(\clX_1, \ttau_1)]^{vir}\times [\clK(\clX_2, \ttau_2)]^{vir})\\
=& a^*l_*[\phi\circ l]([\clK(\clX_1, \ttau_1)]^{vir}\times [\clK(\clX_2, \ttau_2)]^{vir}) \quad (\text{by Proposition \ref{prop:degree_of_l}})\\
=& h_*(\phi\circ l)^!([\clK(\clX_1, \ttau_1)]^{vir}\times [\clK(\clX_2, \ttau_2)]^{vir})\\
=& h_*[\clK(\clX_1\times\clX_2, \ttau)]^{vir}.
\end{split}
\end{equation*}
\end{proo}

\section{Gromov-Witten classes}\label{GW_theory}
In this Section we present consequences of {\bf Theorem \ref{virtual_class_formula}} in Gromov-Witten invariants. 

\subsection{Weighted virtual classes}
We begin with an easy Lemma.
\begin{lemm}\label{product_of_inertia_stack}
Let $\clX_1$ and $\clX_2$ be two Deligne-Mumford stacks. Then there is a natural isomorphism
$$I(\clX_1\times \clX_2)\simeq I(\clX_1)\times I(\clX_2).$$
\end{lemm}
\begin{proo}
This follows easily from the definition of the inertia stacks.
\end{proo}

Gromov-Witten theory of a Deligne-Mumford stack $\clX$ may be defined using the moduli stacks $\clK_{g,n}(\clX, \beta)$. Gromov-Witten classes are defined by intersecting the virtual fundamental class $[\clK_{g,n}(\clX, \beta)]^{vir}$ with cohomology classes of $\clX$ pulled back via the {\em evaluation maps},
$$\overline{ev}: \clK_{g,n}(\clX, \beta)\to \bar{I}(\clX)^{\times n}.$$ A detailed treatment can be found in \cite{AGV06}. Note that the evaluation map $\overline{ev}$ takes values in the {\em rigidified} inertia stack $\bar{I}(\clX)$.

When comparing Gromov-Witten theory of the product stack $\clX_1\times \clX_2$ with Gromov-Witten theory of the factors $\clX_1$ and $\clX_2$, because of {\bf Lemma \ref{product_of_inertia_stack}}, it is more convenient to work in the framework where the evaluation maps take values in the inertia stack. This alternative framework is discussed in \cite{AGVaoqc, AGV06}. We adapt this into our setting, as follows.

Let $\clX$ be a smooth proper Deligne-Mumford stack with projective coarse moduli space $X$. Let $\ttau$ be a gerby $\clX$-graph (see {\bf Definitions \ref{gerby_graph_def}, \ref{gerby_X_graph}}). Over the stack $\clK(\clX, \ttau)$ of $\ttau$-marked twisted stable maps to $\clX$ there are universal marking gerbes $\clG_i, i\in S_\ttau$. The fiber product of all $\clG_i$ over $\clK(\clX, \ttau)$, which we denote by $\clM(\clX, \ttau)$, is the moduli stack of $\ttau$-marked twisted stable maps to $\clX$ {\em with trivialized marked gerbes}. By construction the evaluation map of $\clM(\clX, \ttau)$ fits into the following diagram
\begin{equation*}
\xymatrix{
\clM(\clX, \ttau)\ar[r]^{ev_\ttau}\ar[d]^{\theta_\ttau} & I(\clX)^{\times \#S_{\ttau}}\ar[d]\\
\clK(\clX, \ttau)\ar[r]^{\overline{ev}_\ttau} & \bar{I}(\clX)^{\times \#S_\ttau}. 
}
\end{equation*}
As explained in {\em \cite{AGVaoqc, AGV06}}, in order to define Gromov-Witten theory using $\clM(\clX, \ttau)$ one has to work with the {\em weighted} virtual fundamental class, defined as follows.
\begin{definition}
$$[\clM(\clX, \ttau)]^w:=\left(\prod_{i\in S_\ttau} \gamma(i)\right) \theta_\ttau^*[\clK(\clX, \ttau)]^{vir}.$$
\end{definition}
{\bf Theorem \ref{virtual_class_formula}} has a counterpart for weighted virtual classes. Let $\clX_1~\times~\clX_2$, $\clX_1, \clX_2$, $\ttau$, $\ttau_1$, and $\ttau_2$ be as in Section \ref{virtual_class_result_setup}. Consider the following diagram:
\begin{equation}\label{M_and_K} 
\xymatrix{
\clM(\clX_1\times\clX_2, \ttau)\ar[d]^{\theta'}\ar[dr]^-{\tilde{h}}\ar[drr]\ar@/_1pc/[dd]_{\theta_\ttau} & & \\
\clM\ar[d]^{\theta_{12}}\ar[r]^{h'}\ar@{}[rd]|{\square} & \clP\ar[d]^{\rho}\ar[r]^-{q'}\ar@{}[rd]|{\square} & \clM(\clX_1, \ttau_1)\times \clM(\clX_2, \ttau_2)\ar[d]^{\theta_{\ttau_1}\times \theta_{\ttau_2}}\\  
 \clK(\clX_1\times\clX_2,\ttau)\ar[r]^-{h}& P\ar[r]^-{q}& \clK(\clX_1,\ttau_1)\times \clK(\clX_2,\ttau_2).
 }
\end{equation}
Here the bottom arrows are defined in (\ref{sub_main_diag}). 

\begin{theo}\label{weighted_virtual_class_formula}
\begin{equation}
\tilde{h}_*[\clM(\clX_1\times\clX_2, \ttau)]^{w}=\mathfrak{c}\cdot \Delta^!([\clM(\clX_1, \ttau_1)]^{w}\times [\clM(\clX_2, \ttau_2)]^{w}).
\end{equation}
\end{theo}
\begin{proo}
We calculate 
\begin{equation*}
\begin{split}
\tilde{h}_*[\clM(\clX_1\times \clX_2, \ttau)]^w=&(\prod_{i\in S_\ttau} \gamma(i)) \tilde{h}_*\theta_\ttau^*[\clK(\clX_1\times \clX_2, \ttau)]^{vir}\\
=&(\prod_{i\in S_\ttau} \gamma(i))h'_*\theta'_*\theta'^{*}\theta_{12}^*[\clK(\clX_1\times \clX_2, \ttau)]^{vir}\\
=&(\text{deg}\, \theta')(\prod_{i\in S_\ttau} \gamma(i))h'_*\theta_{12}^*[\clK(\clX_1\times \clX_2, \ttau)]^{vir}\\
=&(\text{deg}\, \theta')(\prod_{i\in S_\ttau} \gamma(i))\rho^* h_*[\clK(\clX_1\times \clX_2, \ttau)]^{vir}.
\end{split}
\end{equation*}
By {\bf Theorem \ref{virtual_class_formula}}, we have 
\begin{equation*}
\begin{split}
&\rho^* h_*[\clK(\clX_1\times \clX_2, \ttau)]^{vir}\\
=&\mathfrak{c} \cdot \rho^* \Delta^!([\clK(\clX_1, \ttau_1)]^{vir}\times [\clK(\clX_2, \ttau_2)]^{vir})\\
=&\mathfrak{c}\cdot \Delta^!(\theta_{\ttau_1}^*[\clK(\clX_1, \ttau_1)]^{vir}\times \theta_{\ttau_2}^*[\clK(\clX_2, \ttau_2)]^{vir})\\
=&\frac{\mathfrak{c}}{\prod_{i\in S_{\ttau_1}}\gamma(i)\prod_{i\in S_{\ttau_2}}\gamma(i)}\Delta^!([\clM(\clX_1, \ttau_1)]^{w}\times [\clM(\clX_2, \ttau_2)]^{w}).
\end{split}
\end{equation*}
We conclude by observing that 
\begin{equation*}
\text{deg}\, \theta'=\frac{\text{deg}\, \theta_\ttau}{\text{deg}\,\theta_{12}}=\frac{\left(\prod_{i\in S_\ttau} \gamma(i)\right)^{-1}}{\left(\prod_{i\in S_{\ttau_1}}\gamma(i)\prod_{i\in S_{\ttau_2}}\gamma(i)\right)^{-1}}.
\end{equation*}
\end{proo}
\subsection{Gromov-Witten classes}
Again let $\clX$ be a proper smooth Deligne-Mumford stack with projective coarse moduli space. Let $\ttau$ be a gerby $\clX$-graph. Denote by $\tau^s$ the absolute stabilization of the underlying modular graph $\tau$. There is a contraction morphism
$$st_{\ttau}: \clM(\clX, \ttau)\to \overline{M}(\tau^s).$$
\begin{definition}
Define a map $$I^\clX_\ttau: H^*(I(\clX), \bbQ)^{\times \# S_{\ttau}}\to H_*(\overline{M}(\tau^s), \bbQ)$$ as follows. For $\omega\in H^*(I(\clX), \bbQ)^{\times \# S_{\ttau}}$, define $$I^\clX_\ttau(\omega):=\left(\prod_{i\in E_\ttau} \gamma(i) \right)st_{\ttau*}(ev_{\ttau}^*(\omega)\cap [\clM(\clX, \ttau)]^w).$$
\end{definition}
We call $I^\clX_\ttau$ the {\em Gromov-Witten class} associated to $\clX$ and $\ttau$.

\subsubsection{Gromov-Witten Cohomological Field Theory}
Let $H$ be a $\bbZ$-graded $\mathbb{Q}$-vector space with a non-degenerate pairing $h$. Following \cite{Mfm}, {\em Chapter 3, Section 4}, a structure of complete (any genus) Cohomological Field Theory (CohTh) on $(H,h)$ is given by  a collection of even linear maps (correlators) 
\begin{eqnarray}\label{correlator_def}
I_{g,n}^V:H^{\otimes n}\to H^*(\overline{M}_{g,n},\mathbb{Q}),
\end{eqnarray}
 defined for any $g\geq 0$ and $n+2g-3\geq 0$, satisfying  the following axioms:
\begin{enumerate}
\item $S_n$-covariance with respect to the natural action od $S_n$ on both sides of  (\ref{correlator_def});

\item Splitting, i.e. compatibility with respect to boundary divisors. Let $\sigma$ be a dual graph with two vertices  of genus $g_ 1$ and $g_2$, defining  a partition $S_1\coprod S_2=\{1,..,n\}$. Let $\phi_\sigma:\overline{M}_{g_1,n_1+1}\times \overline{M}_{g_2,n_2+1}\to \overline{M}_{g,n}$ denote the gluing morphism. Here $n_1=\abs{S_1}$, $n_2=\abs{S_2}$, $g=g_1+g_2$. Then 
\begin{equation*}
\phi_\sigma^*I_{g,n}^V(\gamma_1\otimes ..\otimes \gamma_n)=\epsilon(\sigma)I_{g_1,n_1+1}\otimes I_{g_2,n_2+1}\left(\bigotimes_{p\in S_1}\gamma_p\otimes\Delta\otimes \bigotimes_{q\in S_2}\gamma_q\right),
\end{equation*}
where $\Delta=\sum \Delta_a h^{ab} \otimes \Delta_b$ is the class of the diagonal and $\epsilon(\sigma)$ is the sign of the permutation induced on the odd arguments $\gamma_1$,..,$\gamma_n$;

\item Genus reduction.  For $g\geq 1$, let us denote by $\psi: \overline{M}_{g-1,n+2}\to \overline{M}_{g,n}$ the clutching morphism (\cite{KnudII}, {\em Section 3}). Then we require
\begin{equation*}
 \psi^*I_{g,n}^V(\gamma_1\otimes ..\otimes \gamma_n)=I_{g-1,n+2}(\gamma_1\otimes ..\otimes \gamma_n\otimes\Delta). 
\end{equation*}
\end{enumerate}
\begin{theo}
The classes $I^\clX_{\ttau}$ define a {\em cohomological field theory}.
\end{theo}
\begin{proo}
This follows from compatibility relations of weighted virtual fundamental classes proven in {\em \cite{AGVaoqc}, Proposition 5.3.1}. Consider the splitting axiom of a CohFT. Let $\ttau$ be a gerby $\clX$-graph and $\ttau_1, \ttau_2$ two gerby $\clX$-graphs obtained from $\ttau$ by cutting an edge $\{f, j_\tau f\} \in E_\ttau$ that connects two vertices. Then we can form the following diagram
\begin{equation*}
\xymatrix{
\clM(\clX, \ttau) & \clM(\clX, \ttau_1)\times_{I(\clX)} \clM(\clX, \ttau_2)\ar[l]_-{\pi}\ar[r]\ar[d] &\clM(\clX, \ttau_1)\times \clM(\clX, \ttau_2)\ar[d]\\
& I(\clX)\ar[r]^{\delta} & I(\clX)\times I(\clX).
}
\end{equation*}
Here the fiber product $\clM(\clX, \ttau_1)\times_{I(\clX)} \clM(\clX, \ttau_2)$ is taken using the evaluation maps at flags $f, j_\tau f$. The map $\delta$ is the diagonal map. By definition of the stacks $\clM(\clX, \ttau)$,   $\clM(\clX, \ttau_i)$, the natural map $\pi:  \clM(\clX, \ttau_1)\times_{I(\clX)} \clM(\clX, \ttau_2)\to \clM(\clX, \ttau)$ is the universal gerbe at the node defined by the edge $\{f, j_\tau f\}$, and is of degree $1/\gamma(f)$. The splitting axiom follows from the following consideration:
\begin{equation*}
\begin{split}
&(\prod_{i\in E_\ttau}\gamma(i))[\clM(\clX, \ttau)]^w\\
=& (\prod_{i\in E_\ttau\setminus \{f, j_\tau f\}} \gamma(i)) \gamma(f)^2 \pi_*\pi^*[\clM(\clX, \ttau)]^w\\
=&(\prod_{i\in E_\ttau\setminus \{f, j_\tau f\}} \gamma(i)) \pi_*\delta^!([\clM(\clX, \ttau_1)]^w\times [\clM(\clX, \ttau_2)]^w) \quad \text{by \cite{AGVaoqc}, Proposition 5.3.1}\\
=& \pi_*\delta^!((\prod_{i\in E_{\ttau_1}} \gamma(i))[\clM(\clX, \ttau_1)]^w\times (\prod_{i\in E_{\ttau_2}} \gamma(i))[\clM(\clX, \ttau_2)]^w).
\end{split}
\end{equation*}
The genus-reduction axiom follows from a similar consideration. We omit the details.
\end{proo}
\subsection{Product formula} Again let $\clX_1\times \clX_2$, $\clX_1, \clX_2$, $\ttau$, $\ttau_1$, and $\ttau_2$ be as in Section \ref{virtual_class_result_setup}. We prove here a product formula for Gromov-Witten classes associated to $\clX_1\times \clX_2$, $\clX_1$, and $\clX_2$. 

\begin{theo}\label{GW_class_product_formula}
Let $\omega_1\otimes \omega_2\in H^*(I(\clX_1), \bbQ)^{\times \#S_{\ttau_1}}\otimes H^*(I(\clX_2), \bbQ)^{\times \# S_{\ttau_2}}$. Identify $\omega_1\otimes \omega_2$ with a class in $H^*(I(\clX_1\times \clX_2), \bbQ)^{\times \# S_\ttau}$ via the map $$H^*(I(\clX_1), \bbQ)^{\times \#S_{\ttau_1}}\otimes H^*(I(\clX_2), \bbQ)^{\times \# S_{\ttau_2}}\to H^*(I(\clX_1\times \clX_2), \bbQ)^{\times \# S_\ttau}$$ induced by $F_{\ttau_i}\to F_\ttau$. Then we have in $H_*(\overline{M}(\tau^s), \bbQ)$ $$I^{\clX_1\times \clX_2}_{\ttau}(\omega_1\otimes \omega_2)=I^{\clX_1}_{\ttau_1}(\omega_1)\cup I^{\clX_2}_{\ttau_2}(\omega_2).$$
\end{theo}
\begin{proo}
We refer to diagrams (\ref{main_diag}) and (\ref{M_and_K}). In diagram (\ref{main_diag}) denote the map $P\to \frB^{tw}\to \overline{M}(\tau^s)$ by $\epsilon$. Then $st_\ttau: \clM(\clX_1\times \clX_2, \ttau)\to \overline{M}(\tau^s)$ factors as $st_\ttau= \epsilon\circ \rho\circ \tilde{h}$. Also, the map $st_{\ttau_1}\times st_{\ttau_2}:\clM(\clX_1, \ttau_1)\times \clM(\clX_2, \ttau_2)\to \overline{M}(\tau^s)\times \overline{M}(\tau^s)$ factors as $st_{\ttau_1}\times st_{\ttau_2}=(s\times s)\circ a \circ (\theta_{\ttau_1}\times\theta_{\ttau_2})$. By construction we have $ev_\ttau^*(\omega_1\otimes\omega_2)=\tilde{h}^*q'^*(ev_{\ttau_1}^*(\omega_1)\otimes ev_{\ttau_2}^*(\omega_2))$.  Using these facts and {\bf Theorem \ref{weighted_virtual_class_formula}} we now calculate
\begin{equation*}
\begin{split}
&I^{\clX_1\times \clX_2}_{\ttau}(\omega_1\otimes \omega_2)\\
=&(\prod_{i\in E_\ttau}\gamma(i)) st_{\ttau*}(ev_\ttau^*(\omega_1\otimes \omega_2)\cap [\clM(\clX_1\times \clX_2, \ttau)]^w)\\
=& (\prod_{i\in E_\ttau}\gamma(i)) \epsilon_*\rho_*\tilde{h}_*(\tilde{h}^*q'^*(ev_{\ttau_1}^*(\omega_1)\otimes ev_{\ttau_2}^*(\omega_2))\cap  [\clM(\clX_1\times \clX_2, \ttau)]^w)\\
=&(\prod_{i\in E_\ttau}\gamma(i)) \epsilon_*\rho_*(q'^*(ev_{\ttau_1}^*(\omega_1)\otimes ev_{\ttau_2}^*(\omega_2))\cap \tilde{h}_*[\clM(\clX_1\times \clX_2, \ttau)]^w)\\
=&\mathfrak{c}(\prod_{i\in E_\ttau}\gamma(i)) \epsilon_*\rho_*(q'^*(ev_{\ttau_1}^*(\omega_1)\otimes ev_{\ttau_2}^*(\omega_2))\cap \Delta^!([\clM(\clX_1, \ttau_1)]^{w}\times [\clM(\clX_2, \ttau_2)]^{w}))\\
=&\mathfrak{c}(\prod_{i\in E_\ttau}\gamma(i)) \epsilon_*\rho_*\Delta^!((ev_{\ttau_1}^*(\omega_1)\cap [\clM(\clX_1, \ttau_1)]^{w})\times (ev_{\ttau_2}^*(\omega_2)\cap [\clM(\clX_2, \ttau_2)]^{w}))\\
=&\mathfrak{c}(\prod_{i\in E_\ttau}\gamma(i)) \Delta^*(s\times s)_* a_* (\theta_{\ttau_1}\times\theta_{\ttau_2})_*((ev_{\ttau_1}^*(\omega_1)\cap [\clM(\clX_1, \ttau_1)]^{w})\times (ev_{\ttau_2}^*(\omega_2)\cap [\clM(\clX_2, \ttau_2)]^{w}))\\
=&\mathfrak{c}(\prod_{i\in E_\ttau}\gamma(i)) \Delta^*(st_{\ttau_1*}(ev_{\ttau_1}^*(\omega_1)\cap [\clM(\clX_1, \ttau_1)]^{w})\otimes st_{\ttau_2*}(ev_{\ttau_2}^*(\omega_2)\cap [\clM(\clX_2, \ttau_2)]^{w}))\\
=&I^{\clX_1}_{\ttau_1}(\omega_1)\cup I^{\clX_2}_{\ttau_2}(\omega_2),
\end{split}
\end{equation*}
where in the last line we use {\bf Proposition \ref{prop:degree_of_l}}. The proof is completed. 
\end{proo}

\section{Application to trivial gerbes}\label{trivial_gerbe}
In this Section we consider Gromov-Witten theory of trivial gerbes. Let $\mathcal{X}$ be a proper smooth Deligne-Mumford stack with projective coarse moduli space. Let $G$ be a finite group and $BG$ the classifying stack of $G$. The inertia stack of $BG$ admits the following decomposition 
$$I(BG)=\bigcup_{(g): \text{conjugacy class}} BC_G(g),$$ where $C_G(g)\subset G$ is the centralizer subgroup of $g$ in $G$. This yields a decomposition of the orbifold cohomology groups $$H^*(I(BG), \mathbb{C})=\bigoplus_{(g): \text{conjugacy class}}H^*(BC_G(g), \mathbb{C})=\bigoplus_{(g): \text{conjugacy class}} \mathbb{C}{\bf 1}_{(g)},$$ where ${\bf 1}_{(g)}\in H^0(BC_G(g),\mathbb{C})$ is the standard generator. Let $$I(\mathcal{X})=\bigcup_{i\in \mathcal{I}} \mathcal{X}_i$$ be the decomposition of $I(\mathcal{X})$ into connected components. We know that 
\begin{equation*}
\begin{split}
H^*(I(\mathcal{X}\times BG), \mathbb{C})&\simeq H^*(I(\mathcal{X}), \mathbb{C})\otimes H^*(I(BG), \mathbb{C})\\
&\simeq \bigoplus_{(g):\text{conjugacy class}, i\in \mathcal{I}} H^*(\mathcal{X}_i,\mathbb{C})\otimes{\bf 1}_{(g)}.
\end{split}
\end{equation*}

Let $\beta\in H_2(X, \mathbb{Z})$ be a curve class. Let $\overline{M}_{g,n}(\clX\times BG, \beta)$ be the moduli stack\footnote{This was previously denoted by $\clM(\clX\times BG, \beta)$. In order to be consist with notations in \cite{JK}, we change our notations from now on.} of twisted stable maps with trivialized marked gerbes to $\clX\times BG$ of type $(g, n,\beta)$. Let $(g_1),..., (g_n)$ be conjugacy classes of $G$. Consider integers $k_1,..., k_n\geq 0$ and cohomology classes  $\gamma_j\in H^*(\mathcal{X}_{i_j}, \mathbb{C})$, $1\leq j\leq n$, which we assume to be homogeneous. 
Consider the Gromov-Witten invariant 
$$\langle\prod_{j=1}^n \tau_{k_j}(\gamma_j\otimes {\bf 1}_{(g_j)}) \rangle_{g,n, \beta}^{\mathcal{X}\times BG}:=\int_{[\overline{M}_{g,n}(\mathcal{X}\times BG, \beta)]^{w}}\prod_{j=1}^n ev_j^*(\gamma_j\otimes {\bf 1}_{(g_j)}) \psi_j^{k_j}.$$
It follows from the definition that the descendant classes $\psi_j$ on $\overline{M}_{g,n}(\mathcal{X}\times BG, \beta)$ are obtained as pull-backs of descendant classes on $\overline{M}_{g,n}(\mathcal{X}, \beta)$. This observation together with the property of virtual classes (i.e. {\bf Theorem \ref{weighted_virtual_class_formula}}) imply the following equality:
\begin{equation}
\begin{split}
\langle\prod_{j=1}^n \tau_{k_j}(\gamma_j\otimes {\bf 1}_{(g_j)}) \rangle_{g,n, \beta}^{\mathcal{X}\times BG}=&\langle\prod_{j=1}^n \tau_{k_j}(\gamma_j) \rangle_{g,n, \beta}^{\mathcal{X}}\\
&\times \text{deg}\,(\overline{M}_{g,n}(BG, (g_1),...,(g_n))/\overline{M}_{g,n}).
\end{split}
\end{equation}
Here $\text{deg}\,(\overline{M}_{g,n}(BG, (g_1),...,(g_n))/\overline{M}_{g,n})$ is the degree of the natural map $\overline{M}_{g,n}(BG, (g_1),...,(g_n))\to \overline{M}_{g,n}$. In \cite{JK} this number is denoted by $\Omega_g^G((g_1),...,(g_n))$. According to \cite{JK}, Theorem 3.6 the assignment $$\Lambda_{g,n}^G({\bf 1}_{(g_1)}\otimes...\otimes {\bf 1}_{(g_n)}):=\Omega_g^G((g_1),...,(g_n))$$ defines a CohFT.

Let $\{ V_\alpha\}_{\alpha=1}^r$ be the set of isomorphism classes of irreducible representations of $G$. Denote by $\chi_\alpha$ the character of $V_\alpha$. For $1\leq \alpha\leq r$ define 
$$f_\alpha:=\frac{\text{dim}\, V_\alpha}{|G|}\sum_{(g):\text{conjugacy class}} \chi_\alpha(g^{-1}) {\bf 1}_{(g)},\quad \nu_\alpha:=\left(\frac{\text{dim}\, V_\alpha}{|G|}\right)^2.$$
According to \cite{JK}, Proposition 4.2, we have $\Lambda_{g,n}^G(f_{\alpha_1}\otimes...\otimes f_{\alpha_n})=0$ unless $\alpha_1=...=\alpha_n=:\alpha$, in which case $\Lambda_{g,n}^G(f_\alpha\otimes...\otimes f_\alpha)=\nu_\alpha^{1-g}$. It follows from multi-linearity that 
\begin{equation*}
\begin{split}
&\langle\prod_{j=1}^n \tau_{k_j}(\gamma_j\otimes f_{\alpha_j}) \rangle_{g,n, \beta}^{\mathcal{X}\times BG}\\
=&\begin{cases}
\nu_\alpha^{1-g}\langle\prod_{j=1}^n \tau_{k_j}(\gamma_j) \rangle_{g,n, \beta}^{\mathcal{X}}\quad \text{if }\alpha_1=...=\alpha_n=:\alpha,\\
0, \quad \text{otherwise}.
\end{cases}
\end{split}
\end{equation*}

This may be formulated in terms of generating functions. Let $\{ \gamma_j\}_{1\leq j\leq \text{rank}\,H^*(I(\mathcal{X}),\mathbb{C})} \subset H^*(I(\mathcal{X}),\mathbb{C})$ be an additive basis consisting of homogeneous elements. Let $\{t_{j, k}| 1\leq j\leq \text{rank}\,H^*(I(\mathcal{X}),\mathbb{C}), k\geq 0\}$ and $\{t_{\alpha, j,k}|1\leq \alpha\leq r, 1\leq j\leq \text{rank}\,H^*(I(\mathcal{X}),\mathbb{C}), k\geq 0\}$ be two sets of variables. We may consider generating functions of genus $g$ descendant Gromov-Witten invariants:
\begin{equation*}
\begin{split}
&\mathcal{F}_{\mathcal{X}\times BG}^g(\{t_{\alpha, j,k}\}_{1\leq \alpha\leq r, 1\leq j\leq \text{rank}\,H^*(I(\mathcal{X}),\mathbb{C}), k\geq 0}; Q)\\
:=&\sum_{\overset{n, \beta}{\alpha_1,...,\alpha_n; j_1,...,j_n; k_1,...,k_n}}\frac{Q^\beta}{n!}\prod_{l=1}^n t_{\alpha_l, j_l, k_l}\langle\prod_{l=1}^n\tau_{k_l}(\gamma_{j_l}\otimes f_{\alpha_l})\rangle_{g,n,\beta}^{\mathcal{X}\times BG};\\
&\mathcal{F}_{\mathcal{X}}^g(\{t_{j, k}\}_{1\leq j\leq \text{rank}\,H^*(I(\mathcal{X}),\mathbb{C}), k\geq 0}; Q)\\
:=&\sum_{\overset{n, \beta}{j_1,...,j_n; k_1,...,k_n}}\frac{Q^\beta}{n!}\prod_{l=1}^n t_{ j_l, k_l}\langle\prod_{l=1}^n\tau_{k_l}(\gamma_{j_l})\rangle_{g,n,\beta}^{\mathcal{X}}.\\
\end{split}
\end{equation*}

We thus have 
\begin{equation}
\begin{split}
&\mathcal{F}_{\mathcal{X}\times BG}^g(\{t_{\alpha, j,k}\}_{1\leq \alpha\leq r, 1\leq j\leq \text{rank}\,H^*(I(\mathcal{X}),\mathbb{C}), k\geq 0}; Q)\\
=&\sum_{\alpha=1}^r\nu_\alpha^{1-g}\mathcal{F}_{\mathcal{X}}^g(\{t_{\alpha, j, k}\}_{1\leq j\leq \text{rank}\,H^*(I(\mathcal{X}),\mathbb{C}), k\geq 0}; Q).
\end{split}
\end{equation}
In terms of the {\em total descendant potential} $\mathcal{D}:=\exp(\sum_{g\geq 0} \hbar^{g-1}\mathcal{F}^g)$, we have 
\begin{equation}
\begin{split}
&\mathcal{D}_{\mathcal{X}\times BG}^g(\{t_{\alpha, j,k}\}_{1\leq \alpha\leq r, 1\leq j\leq \text{rank}\,H^*(I(\mathcal{X}),\mathbb{C}), k\geq 0}; Q)\\
=&\prod_{\alpha=1}^r\mathcal{D}_{\mathcal{X}}^g(\{t_{\alpha, j, k}\}_{1\leq j\leq \text{rank}\,H^*(I(\mathcal{X}),\mathbb{C}), k\geq 0}; Q)_{\hbar\mapsto \hbar/\nu_\alpha}.
\end{split}
\end{equation}
This verifies the decomposition conjecture for $\mathcal{X}\times BG$. 


\bibliographystyle{alpha}
\bibliography{el_biblio.bbl}

\end{document}